\documentclass[11pt]{article}
\usepackage[margin=1in]{geometry}  
\usepackage{graphicx}
\usepackage{amssymb}             
\usepackage{amsmath}
\usepackage{slashed}            
\usepackage{amsfonts}              
\usepackage{amsthm}
\usepackage{breqn}
\usepackage{verbatim}              
\usepackage{enumitem}
\usepackage{fancyhdr}
\theoremstyle{definition}
\newtheorem{thm}{Theorem}[section]

\newtheorem{mydef}{Definition}[section]
\newtheorem{lem}[thm]{Lemma}
\newtheorem{prop}[thm]{Proposition}

\newtheorem*{rmk}{Remark}

\newcommand{\RR}{\mathbb{R}}      

\newcommand{\q}{\quad}
\newcommand{\p}{\partial}
\newcommand{\DD}{\mathcal{D}}
\newcommand{\curl}{\text{curl}\,}
\newcommand{\lei}{L^{2}(\mathcal{D}_{t})}
\newcommand{\leb}{L^{2}(\p\mathcal{D}_{t})}
\newcommand{\lli}{L^{2}(\Omega)}
\newcommand{\llb}{L^{2}(\p\Omega)}
\newcommand{\den}{{\bar\rho}_{0}}
\newcommand{\Dt}{{\p}_{t}+v^{k}{\p}_{k}}

\newcommand{\nab}{\nabla}

\newcommand{\lap}{\Delta}
\newcommand{\di}{\text{div}\,}
\newcommand{\detg}{\text{det}\,g}
\newcommand{\cnab}{\overline{\nab}}
\newcommand{\cp}{\overline{\partial}}
\newcommand{\dx}{\,dx}
\newcommand{\vol}{\text{vol}\,}
\newcommand{\lee}{\langle\langle}
\newcommand{\ree}{\rangle\rangle}
\newcommand{\symdot}{\tilde{\cdot}}
\newcommand{\linf}{L^{\infty}(\p\Omega)}
\newcommand{\kk}{\kappa}
\newcommand{\E}{\widehat{E}}
\numberwithin{equation}{section}

\begin{document}
\title{A priori Estimates for the Compressible Euler Equations for a Liquid with Free Surface Boundary and the Incompressible Limit}
\author{Hans Lindblad and Chenyun Luo}
\date{}
\maketitle
\begin{abstract}
 In this paper, we prove a new type of energy estimates for the compressible Euler's equation with free boundary, with a boundary part and an interior part. These can be thought of as a generalization of the energies in Christodoulou and Lindblad \cite{CL} to the compressible case and do not require the fluid to be irrotational. In addition, we show that our estimates are in fact uniform in the sound speed $\kk$. As a consequence, we obtain convergence of solutions of compressible Euler equations with a free boundary to solutions of the incompressible equations,
  generalizing the result of Ebin \cite{Eb} to when you have a free boundary. In the incompressible case our energies reduces to those in \cite{CL} and our proof in particular gives a simplified proof of the estimates in \cite{CL} with improved error estimates. Since for an incompressible irrotational liquid with free surface there are small data global existence results our result  leaves open the possibility of long time existence also for slightly compressible liquids
  with a free surface.
\end{abstract}
\tableofcontents

\section{Introduction}
\quad We consider Euler equations
\begin{align}
\begin{cases}
(\Dt) v = -{\rho}^{-1}\p p \q \text{in}\, \mathcal{D},\\ \label{E}
(\Dt) \rho+\rho \, \di v =0 \q \text{in}\, \mathcal{D}.
\end{cases}
\end{align}
describing the motion of a perfect compressible fluid in vacuum, where $\p_{i} = \frac{\p}{\p x_i}$, $\di v=\p_kv^k$, and $v=(v_1,\cdots,v_n)$ and $\mathcal{D} = \cup_{0\leq t\leq T}\{t\}\times \mathcal{D}_t$, $\DD_t\subset \RR^n$, and the density $\rho$ are to be determined. Here, $v^{k} = \delta^{ij}v_j = v_k$, and we have used the summation convention on repeated upper and lower indices. The pressure $p=p(\rho)$ is assumed to be a given strictly increasing smooth function of the density. The boundary $\p \DD_t$ moves with the velocity of the fluid particles at the boundary. The fluid moves in the vacuum so the pressure $p$ vanishes in the exterior and hence on the boundary. We therefore requires the boundary condition on $\p \DD = \cup_{0\leq t\leq T}\{t\} \times \p\DD_t$ to be
\begin{align}
\begin{cases}
(\Dt) |_{\p\DD} \in T(\p\DD),\\ \label{BC}
p=0 \q \text{on}\,\p\DD.\\
\end{cases}
\end{align}
Since the pressure is constant on the boundary the density has to be a constant $\den\geq0$ on the boundary. We assume in fact $\den>0$, which is in the case of a liquid (as apposed to a gas). Hence,
\begin{align}
p(\den)=0, \q p'(\rho)>0, \q \text{for}\q \rho\geq\den, \label{pbdy}
\end{align}
where, for the sake of simplicity, we further assume $\den =1$. \\
\indent Given a bounded domain $\DD_0\subset\RR^n$, that is homeomorphic to the unit ball, and initial data $v_0$ and $\rho_0$, we want to find a set $\DD$, a vector field $v$ and a function $\rho$, solving (\ref{E})-(\ref{BC}) and satisfying the initial conditions
\begin{align}
\begin{cases}
\{x:(0,x)\in\DD\} = \DD_0,\\
v=v_0,\rho=\rho_0 \q \text{on}\,\{0\}\times\DD_0. \label{IC}
\end{cases}
\end{align}

Since the pressure is an increasing function of the density one can alternatively think of the density as a function of the pressure. By thinking of the density as a function of the pressure (or rather the enthalphy, see below) the incompressible case can be thought of as a special case of constant density function. This point of view that we will take is needed in order to be able to pass to the incompressible limit for the free boundary problem.

\subsection{Enthalpy form}
\q Let $D_t = \Dt$ be the material derivative. We introduce the enthalpy $h$ to be the strictly increasing function of the density; $h(\rho)=\int_{1}^{\rho}p'(\lambda)\lambda^{-1}\,d\lambda$. Since the enthalphy is a strictly increasing function of the density we can alternatively think of the density as a function of the enthalphy $\rho(h)$:
\begin{equation}
\rho(0)=\den,\qquad \rho^\prime(h)>0,\quad \text{for}\quad h\geq 0.
\end{equation}
We define $e(h)=\log \rho(h)$. Under these new variables, (\ref{E})-(\ref{IC}) can be re-expressed as
\begin{align}
\begin{cases}
D_t v = -\p h\q \text{in}\, \mathcal{D},\\
\di v = -D_t e(h)=-e'(h)D_t h\q \text{in}\, \mathcal{D}. \label{EE}
\end{cases}
\end{align}
Together with initial and boundary conditions
\begin{align}
\begin{cases}
\{x:(0,x)\in\DD\} = \DD_0,\\
v=v_0,h=h_0 \q \text{on}\,\{0\}\times\DD_0, \label{EIBC}
\end{cases}
\begin{cases}
D_t|_{\p \DD}\in T(\p\DD),\\
h = 0 \q \text{on}\,\p\DD,
\end{cases}
\end{align}
(\ref{EE}) looks exactly like the incompressible Euler's equations, where $h$ takes the position of $p$ and $\di v$ is no longer $0$ but determined by $h$. On the other hand, we would like to impose the following natural conditions on $e(h)$ for some fixed constant $c_0$
\begin{equation}\label{econd}
|e^{(k)}(h)|\leq c_0,\qquad\text{and}\qquad
|e^{(k)}(h)|\leq c_0\sqrt{e'(h)}, \qquad\text{for $k\leq 6$}.
\end{equation}

\indent In order for the initial boundary problem (\ref{EE})-(\ref{EIBC}) to be solvable the initial data has to satisfy certain compatibility conditions at the boundary. But the second equation in (\ref{E}),(\ref{pbdy}) implies that $\di v|_{\p\DD}=0$. We must therefore have $h_0|_{\p\DD_{0}}=0$ and $\di v_0|_{\p\DD_0}=0$, which is the zero-th compatibility condition. Furthermore, $m$-th order compatibility condition can be expressed as
\begin{align}
(\Dt)^j h|_{\{0\}\times\p\DD_{0}}=0,\q j = 0,\cdots,m. \label{cpt cond}
\end{align}
\indent Let $N$ be the exterior unit normal to the free surface $\p\DD_t$. We will prove a priori bounds for (\ref{EE})-(\ref{EIBC}) in Sobolev spaces under the assumption
\begin{align}
\nab_{N}h\leq -\epsilon<0,\q \text{on}\,\p\DD_t, \label{RE}
\end{align}
where $\nab_{N}=N^i\p_i$ and $\epsilon>0$ is a constant. (\ref{RE}) is a natural physical condition. It says that the pressure and hence the enthalpy is larger in the interior than at the boundary. The system (\ref{EE})-(\ref{EIBC}) is ill-posed in absence of (\ref{RE}), an easy counter-example can be found in \cite{CL}.

\subsection{History and background}
\indent Euler equations involving free-boundary has been studied intensively by many authors. The first break through in solving the well-posedness for the incompressible and irrotational problem for general data came in the work of Wu \cite{W1,W2} who solved the problem in both two and three dimensions. For the general incompressible problem with nonvanishing curl Christodoulou and Lindblad \cite{CL} were the first to obtain the energy estimates assuming the Taylor sign condition. For the compressible problem, Lindblad \cite{L} later proveed local well-posedness for the general problem modelling the motion of a liquid via Nash-Moser iteration. On the other hand, Coutand-Lindblad-Shkoller \cite{CLS} and Jang-Masmoudi \cite{JM} obtained the energy estimates and well-posedness for the general problem modelling the motion of gas. It is worth mentioning that D. Ebin \cite{Eb}, and Ebin-Disconzi \cite{DE} proved the solutions of the compressible equations converges to the solutions of the incompressible equation in Sobolev norms as the sound speed goes to infinity, but within a domain with fixed boundary. But no previous incompressible limit result involving free boundary is known. Our result allows us to approximate slightly compressible liquid by the incompressible liquid in both $2$D and $3$D, for which global (in time) solution is known to exist (e.g. \cite{GMS, HIT,IT,IP,W3,W4}).

\indent In this paper, we generalize the method used by Chistodoulou and Lindblad \cite{CL}. In our proof, $\curl v$ appears to be of lower orders. In addition, our method is regardless of spatial dimensions. The energy constructed in this paper contains interior and boundary parts, where the interior part controls the velocity and the enthalpy in Sobolev norms. The boundary part contains projected spatial derivatives, which controls the second fundamental form of the moving boundary. The use of projected derivatives on the boundary is crucial due to the loss of regularity when estimating on the boundary, i.e., the trace theorem \cite{Ev}, and the use of the tangential part of derivatives on the boundary compensates the loss.

\subsection{Energy conservation and higher order energies}
\indent The boundary condition $p|_{\p\Omega}=0$ leads to that the zero-th order energy is conserved, i.e., let
\begin{align}
E_0(t)= \frac{1}{2}\int_{\DD_t}\rho|v|^2\,dx +\int_{\DD_t}\rho \,Q(\rho)\,dx,
\end{align}
where $Q(\rho)=\int_{1}^{\rho}p(\lambda)\lambda^{-2}\,d\lambda $.
A direct computation yields
\begin{multline}
\frac{d}{dt}E_0(t)=\int_{\DD_t}\rho D_t v \cdot v\dx+\int_{\DD_t}p(\rho)D_t\rho \rho^{-1}\dx
=-\int_{\DD_t}\p p\cdot v\dx+\int_{\DD_t}p(\rho)D_t\rho \rho^{-1}\dx\\
=-\int_{\DD_t}p\,\di v\dx+\int_{\DD_t}p(\rho)D_t\rho \rho^{-1}\dx=0.\label{energyconservation}
\end{multline}

 In order too define  higher order energies we introduce a positive definite quadratic form $Q$ on $(0,r)$ tensors, which, when restricted to the boundary, is the inner product of the tangential components, i.e.,
 \begin{equation}
 Q(\alpha,\beta)=\langle \Pi \alpha,\Pi\beta\rangle,\qquad\text{on } \p\DD_t,
 \end{equation}
where the projection of a  $(0,r)$ tensor to the boundary is defined by
$$
(\Pi \alpha)_{i_1,\cdots,i_r}=\gamma_{i_1}^{j_1}\cdots\gamma_{i_r}^{j_r}\alpha_{j_1,\cdots,j_r},
\qquad\text{where } \gamma_i^j=\delta_i^j-\mathcal{N}_i \mathcal{N}^j,
$$
and $\mathcal{N}$ is the unit normal to $\p\DD_t$. To be more specific, in the interior we define
\begin{align}
Q(\alpha,\beta)= q^{i_1j_1}\cdots q^{i_rj_r}\alpha_{i_1\cdots i_r}\beta_{j_1\cdots j_r},
\end{align}
where
\begin{align*}
q^{ij}=\delta^{ij}-\eta(d)^2\mathcal{N}^i\mathcal{N}^j,\\
d(x) = dist(x,\p\DD_t),\\
\mathcal{N}^i=-\delta^{ij}\p_j d.
\end{align*}
Here $\eta$ is a smooth cut-off function satisfying $0\leq \eta(d)\leq1$, $\eta(d)=1$ when $d\leq\frac{d_0}{4}$ and $\eta(d)=0$ when $d>\frac{d_0}{2}$. $d_0$ is a fixed number that is smaller than the injective radius of the normal exponential map $l_0$,defined to be the largest number $l_0$ such that the map
\begin{align}
\p\DD_{t}\times(-l_0,l_0)\to\{x:dist(x,\p\DD_t)<l_0\}, \label{inj rad}
\end{align}
given by
\begin{align}
(\bar{x},l)\to x=\bar{x}+l\mathcal{N}(\bar{x}),
\end{align}
is an injection.\\
\indent The  higher order energies we propose are
\begin{equation}
E_r = \sum_{s+k=r} E_{s,k}+K_r+W_{r+1}^2\label{Er},\qquad E_r^*=\sum_{r'\leq r} E_{r'},
\end{equation}
 where
\begin{multline}
E_{s,k}(t)=\frac{1}{2}\int_{\DD_t}\rho \delta^{ij}Q(\p^sD_t^kv_i,\p^sD_t^{k}v_j)\dx+\frac{1}{2}\int_{\DD_t}\rho e'(h)Q(\p^sD_t^k h,\p^s D_t^kh)\dx\\
+\frac{1}{2}\int_{\p\DD_t}\rho Q(\p^sD_t^k h,\p^sD_t^k h)\nu\,dS,  \label{Esk}
\end{multline}
where $\nu = (-\nab_{N}h)^{-1}$ and
\begin{align}
K_r(t) &= \int_{\DD_t}\rho|\p^{r-1}\curl v|^2\dx \label{K_r},\\
W_r(t) &= \frac{1}{2}||\sqrt{e'(h)}D_t^rh||_{\lei}+\frac{1}{2}||\nab D_t^{r-1}h||_{\lei} .\label{Er time}
\end{align}

Here $W_r$ is the (higher order) energy for the wave equation
\begin{align}
D_t^2e(h) -\lap h=(\p_i v^j)(\p_j v^i),\label{waveequation}
\end{align}
which is obtained by commuting divergence through the first equation of (\ref{E}) using
\begin{equation}
[D_t, \p_i] = -(\p_iv^j)\p_j.\label{lowestcommutator}
\end{equation}
Similarly it follows that we have a transport equation for the curl
\begin{align}
D_t\curl_{ij}v=-(\p_iv^k)(\curl_{kj}v)+(\p_jv^k)(\curl_{ki}v).
\end{align}

Although the energies $E_r$ only control the tangential components, the fact that
we also control the divergence $W_{r+1}^2$ (through $\di v=-D_t e(h)$) and the curl $K_r$ allows us to control all components. In fact, by a Hodge type decomposition
\begin{align}
|\p v| \lesssim |\overline{\p}v|+|\di v|+|\curl v|,
\end{align}
where the tangential derivatives are given by $\overline{\p}h=\Pi\p h$.

The boundary term in \eqref{Esk} and  $\nu$ are constructed to exactly cancel a boundary term coming from integration by parts in the interior as in \eqref{energyconservation}, as will be explained in section  1.5. Moreover the projection in the boundary term is needed to make it lower order in  space derivatives of $h$.
In fact, since $h$ vanishes on the boundary so does the tangential derivatives $\overline{\p}h=\Pi\p h$ and similarly $\Pi \p^r h=O(\p^{r-1}h)$ is lower order.

Moreover if $|\nabla_{\mathcal{N}} h|\geq \epsilon>0$ then the boundary term gives an estimate
for the regularity of the boundary. In fact, one can show that if $q$ vanishes on the boundary then
\begin{equation}
\Pi \p^r q = (\overline{\p}^{r-2}\theta)\nab_{\mathcal{N}} q+O(\p^{r-1}q)+O(\overline{\p}^{r-3}\theta),\label{projest0}
\end{equation}
where $\theta$  is the second fundamental form of the boundary and $\overline{\p}$ stand for tangential derivatives, so
\begin{equation}
\| \overline{\p}^{r-2}\theta\|_{\leb}^2 \leq \frac{C}{\epsilon} E_r^*
+ C \sum_{r'\leq r-1} \| \p^{r'} h\|_{\leb}^2.\label{projest}
\end{equation}

Because of the bound on the second fundamental form energies in fact control all components
\begin{align*}
||v||_{r,0} := \sum_{k+s=r, k<r} ||\p^s D_t^k v||_{\lei},\\
||h||_r := \sum_{k+s=r, k<r}||\p^s D_t^k h||_{\lei} + ||\sqrt{e'(h)}D_t^rh||_{\lei},\\
\lee h\ree_r :=\sum_{k+s=r} ||\p^s D_t^kh||_{\leb},
\end{align*}
in the interior and on the boundary.  Using elliptic estimates (see section 3) one can show that
\begin{align}
||v||_{r,0}^2+||h||_r^2 \leq C(K,M,c_0,\vol \DD_t)E_r^* \label{intro int  est},\\
||D_th||_r^2+\lee h\ree_r^2 \leq C(K,M,c_0,\frac{1}{\epsilon}, \vol \DD_t, E_{r-1}^*)E_r^* \label{intro bdy est}.
\end{align}

\subsection{The main results}
We prove energy estimates implying that the higher order energies remain bounded as long as certain a priori assumptions are true.
\thm
 Let $(v, h)$ be the solution for  (\ref{EE})-(\ref{EIBC}) and $E_r$ be defined as in \eqref{Er}. Then there are continuous functions $C_r$ such that
\begin{align}
\Big|\frac{dE_r(t)}{dt}\Big|\leq C_r(K,\frac{1}{\epsilon},M,c_0, \vol\DD_t,E_{r-1}^*)E_r^*(t),\label{higherorderenergyestimate}
\end{align}
if $0\leq r\leq 4$ provided that the assumptions \eqref{econd} on $e(h)$ hold and
\begin{align}
|\theta|+\frac{1}{l_0}\leq K,\q\text{on}\,\p\DD_t,\label{geometry_bound}\\
-\nab_{N}h\geq \epsilon>0,\q\text{on}\,\p\DD_t,\label{RT sign}\\
1\leq |\rho|\leq M,\q\text{in}\,\DD_t, \label{bootstrap rho}\\
|\p v|+|\p h|+|\p^2 h|+|\p D_th|\leq M,\q\text{in}\,\DD_t. \label{v,h}\\
|D_t h|+|D_t^2 h|\leq M,\q\text{in}\,\DD_t. \label{Dtth}
\end{align}
The bounds (\ref{geometry_bound}) gives us control of geometry of the free surface $\p\DD_t$. A bound for the second fundamental form $\theta$ gives a bound for the curvature of $\p\DD_t$, and a lower bound for the injectivity radius of the exponential map $l_0$ measures how far off the surface is from self-intersecting. Note that for the compressible Euler's equations the bounds (\ref{bootstrap rho})-(\ref{Dtth}) together with the second equation of (\ref{EE}) and \eqref{waveequation} imply the bounds \eqref{Dtth}. We only include these bounds here because we need them to hold uniformly to pass to the incompressible limit.
It follows from \eqref{higherorderenergyestimate} that the energies $E_r(t)$ are bounded as long as the apriori $L^\infty$ bounds above hold. On the other hand it follows from the energy bounds if $r\geq 4$ and $n\leq 3$ that the
a priori $L^\infty$ bounds hold up to some some small positive time $t\leq T$ (depending only on the initial energy and $L^\infty$ bounds) if slightly stronger bounds hold initially (Proposition \ref{a priori estimate strong}).

The above energy bounds remain valid uniformly as the sound speed goes to infinity (Theorem \ref{uniform energybounds} and Proposition \ref{uniform kk}).
The sound speed $\kk$ is defined be viewing $\{p_{\kk}(\rho)\}$ as a family parametrized by $\kk\in\RR^+$, such that for each $\kk$ we have
$$
p_{\kk}'(\rho)|_{\rho=1}=\kk.
$$
Under this setting, we consider the Euler equations
\begin{align}
\begin{cases}
D_t v_{\kk} = -\p h_{\kk}\q \text{in}\, \mathcal{D},\\
\di v_{\kk} = -D_t e_{\kk}(h_\kk)\q \text{in}\, \mathcal{D} \label{EEkk}.
\end{cases}
\end{align}
We view the density as a function of the enthalpy, i.e., $\rho_{\kk} = \rho_{\kk}(h)$. We further assume that $\rho_{\kk}(h)$ and $e_{\kk}(h):=\log \rho_{\kk}(h)$ satisfies
\begin{equation}
\rho_{\kk}(h)\to 1, \quad\text{and}\quad e_{\kk}(h)\to 0,\qquad \text{as}\quad
\kk\to\infty
\end{equation}
and for some fixed constant $c_0$
\begin{equation}\label{econdkappa}
|e^{(k)}_{\kk}(h)|\leq c_0,\qquad\text{and}\qquad
|e^{(k)}_{\kk}(h)|\leq c_0\sqrt{e^{\,\prime}_{\kk}(h)}, \qquad\text{for $k\leq 6$}.
\end{equation}
\thm \label{intro convergence} Let $u_0$ be a divergence free vector field such that its corresponding pressure $p_0$, defined by $\lap p_0 = -(\p_i u_0^k)(\p_k u_0^i)$ and $p_0\big|_{\p\DD_0}=0$, satisfies the physical condition $-\nab_N p_0\big|_{\p\DD_0} \geq \epsilon >0$. Let $(u, p)$ be the solution of the incompressible free boundary Euler equations with data $u_0$, i.e.
$$
\rho_0 D_t u = -\p p,\q \di u=0,\qquad p|_{\p\DD_0}=0,\quad u|_{t=0}=u_0
$$
with the constant density $\rho_0 =1$. Furthermore, let $(v_{\kk}, h_{\kk})$ be the solution for the compressible Euler equations \eqref{EEkk}, with the density function $\rho_{\kk} :h\rightarrow \rho_{\kk}(h)$, and the initial data $v_{0\kk}$ and $h_{\kk}|_{t=0} = h_{0\kk}$, satisfying the compatibility condition \eqref{cpt cond} up to order $r+1$, as well as the physical sign condition \eqref{RE}. Suppose that $\rho_{\kk} \to \rho_0=1$, $v_{0\kk}\to u_0$ and $h_{0\kk}\to p_0$ as $\kk\to\infty$, such that $E_{r,\kk}^*(0)$ is bounded uniformly independent of $\kk$, then
$$
(v_{\kk},h_{\kk})\to (u,p).
$$

The uniform (with respect to the sound speed) a priori bounds are due to that our estimates do not depend on the lower bound of $e_{\kk}^{(k)}(h)$, which goes to $0$ as $\kk\to\infty$. In addition, apart from the coefficient in front of the highest order time derivative our energy does not depend in crucial way on $\kk$ but uniformly (as $\kk\to \infty$) control the corresponding norms of all but the highest order time derivative. This leads to that the a priori
$L^\infty$ bounds also hold uniformly and the norms are bounded uniformly up to a fixed time. The convergence of solutions for the compressible Euler equations to the solution for the incompressible equations then follows.  Furthermore, our energy estimate can be slightly modified so that it can be carried over to the incompressible Euler equations. We refer the remarks after Proposition \ref{uniform kk} for details.

Finally, in section 7 we show that for every divergence free $u_0\in H^{s}, s\geq 5$, there are initial data $v_{0\kk}$, $h_{0\kk}$ for the compressible equations for large $\kk$ satisfying the required number of compatibility conditions and converging in our energy norm to the incompressible data as $\kk\to\infty$, so that the assumptions in the previous theorem hold, and hence the incompressible limit exist. We show

\thm \label{convergence of data} Let $u_0$ and $p_0$ are the initial data for the incompressible Euler equations defined in Theorem \ref{intro convergence}, and we further assume $u_0\in H^s, s\geq 5$. Let $\rho_\kappa(h)\sim\rho_0+h/\kappa$. Then there exists initial data $v_{0\kk}$ and $h_{0\kk}$ satisfying the compatibility condition \eqref{cpt cond} up to order $r+1$, such that $v_{0\kk}\to u_0$,  $h_{0\kk}\to p_0$ as $\kk\to \infty$, and $E_{r,\kk}^*(0)$ is uniformly bounded for all $\kk$.

In addition, given the compressible initial data in terms of the enthalpy $h$ seems necessary for the case when there is a moving domain. This is in order that the compaibility conditions should be satisfied. Also, otherwise if the data is in terms of the density $\rho$, there would be a limited choice of the initial velocity $v_0$, and apart from this, the Taylor sign condition \eqref{RE} may fail on the moving boundary. We refer the remarks after Theorem \ref{existence of data in bounded initial domain} for details. Above we make a particular simple choice of $\rho_\kappa(h)$ just to make the calcualtions simpler.

\rmk
 For simplicity we shall only prove \eqref{higherorderenergyestimate} for $r\leq 4$, which will be sufficient to get back the apriori bounds and close the argument for $n=2,3$. Our method can however be used to prove the energy bound for all order $r$. Since existence in $H^N$ for some large $N$ was shown in \cite{L} using Nash-Moser iteration, one should expect that the solution $v_{\kk}, h_{\kk}$ exist in some fixed time interval $[0,T]$ as long as the energy bounds of order $r\geq N$ hold.

\rmk
One could try to define $E_r$ with $\widehat{D}_t:= \sqrt{e_{\kk}'(h)}D_t$ in place of $D_t$, as an analog to Disconzi-Ebin \cite{DE} and Ebin \cite{Eb}. However, the resulting energies are too weak to control the evolution of
$$
\mathcal{E}(t):=|(\nab_N h(t,\cdot))^{-1}|_{\linf}.
$$
This is due to the a priori assumption $|\p D_th| \leq M$ has to be replaced by $|\p \widehat{D}_th|\leq M$, which is weaker since $e_{\kk}'(h) \to 0$. Although our energies are
stronger we show that for every incompressible data there are data for the compressible equations converging in our energy norm.

\rmk We can alternatively use the modified energies defined as
\begin{equation}
\hat{E}_r = \sum_{s+k=r,\,k\leq r-2} E_{s,k}+K_r+W_{r+1}^2,
\end{equation}
which reduces the number of time derivatives involved in $E_r$. The statement of the main theorem still holds with $E_r$ replaced by $\E_r$.

\subsection{Outline of the proof of the higher order energy estimate \eqref{higherorderenergyestimate}}
We conclude the introduction by showing how the time derivative of the interior terms of the energy to leading order cancel each other after integrating by parts modulo a boundary term that in turn is  to leading order canceled by the time  derivative of the  boundary term.
Let $s+k=r$, we have
\begin{multline}
\frac{d}{dt}E_{s,k} = \int_{\DD_t}\rho \delta^{ij}Q(\p^sD_t^kv_i,D_t \p^sD_t^{k}v_j)\dx
+\int_{\DD_t}\rho e'(h)Q(\p^sD_t^k h,D_t \p^s D_t^kh))\dx\\
+\int_{\p\DD_t}\rho Q(\p^sD_t^k h,D_t \p^sD_t^k h)\nu \,dS+\dots,,
\end{multline}
where  the dots stand for lower order terms.
Using the commutator $[D_t, \p_i] = -(\p_iv^j)\p_j$ and the equation $D_t v_i = -\p_i h$ we get
\begin{align}
D_t\p^sD_t^k v_j &= -\p^sD_t^k\p_j h +\dots=-\p_j \p^sD_t^k h +\dots,\\
D_t\p^rh + (\p_jh)\p^rv^j &= \p^rD_th+\dots,\label{boundarycommutator}\\
D_t\p^sD_t^k h &= \p^sD_t^{k+1}h+\dots.
\end{align}
Hence,
\begin{multline}
\frac{d}{dt}E_{k,s}= -\int_{\DD_t}\rho (\delta^{ij}Q(\p^sD_t^kv_i,\p_j\p^{s}D_t^{k}h)\dx+\int_{\DD_t}\rho e'(h)Q(\p^sD_t^k h,\p^{s} D_t^{k+1} h)\dx\\
+\int_{\p\DD_t}\rho Q(\p^sD_t^k h,D_t\p^sD_t^k h)\nu\,dS
+\dots.
\end{multline}
Now, if we integrate by part in the first term, we get
\begin{multline}
\frac{d}{dt}E_{k,s} =\int_{\DD_t}\rho \delta^{ij}Q(\p_j \p^sD_t^k v_i,\p^sD_t^kh)\dx+\int_{\DD_t}\rho e'(h)Q(\p^sD_t^k h,\p^{s} D_t^{k+1} h)\dx\\
+\int_{\p\DD_t}\rho Q(\p^sD_t^k h,D_t\p^sD_t^k h-\nu^{-1}N_i\p^sD_t^kv^i)\nu\,dS
+\dots.
\end{multline}
The terms in the first line cancel each other (up to lower-order terms) since
$\delta^{ij}\p_j \p^sD_t^k v_i=\p^sD_t^k\delta^{ij}\p_j v_i+\dots$ and $\di v = -e'(h) D_th$.

Because our total energy of order $r$ contain estimates of  more time derivatives than space derivatives the most problematic case in which we need to estimate the boundary term above is
when $s=r$ and $k=0$. Using \eqref{boundarycommutator} we see hence see that we are left with
\begin{equation}
\frac{d}{dt}E_{r,0}
=\int_{\p\DD_t}\rho Q(\p^r h,\p^r D_t h -\p_i h\, \p^r v^i -\nu^{-1}N_i\p^r v^i)\nu\,dS
+\dots.
\end{equation}
We have choose $\nu$ to exactly cancel the leading order term at the boundary in this case.
Since $-\nu^{-1}N_i = \p_ih$, the first term on the second line is inner product of $||\Pi\p^r h||_{\leb}$ and
plus the sum of the inner products of $||\Pi\p^r D_t h||_{\leb}$, which due to \eqref{projest0}-\eqref{projest}
we are able to control.

The proof of the energy estimate for Euler's equations outlined above is given in section 5. The proof of the energy estimate for the wave equation in section 4 and the elliptic bounds in section 3.

\section{Lagrangian coordinate, covariant differentiation and metric, regularity of the boundary }
\indent Let us first introduce Lagrangian coordinate, under which the boundary becomes fixed. Let $\Omega$ be the unit ball in $\RR^n$, and let $f_0:\Omega\to\DD_0$ to be a diffeomorphism. The Lagrangian coordinate $(t,y)$ where $x=x(t,y)=f_t(y)$ are given by solving
\begin{align}
\frac{dx}{dt}=v(t,x(t,y)), \q x(0,y)=f_0(y),\q y\in\Omega. \label{change}
\end{align}
The boundary becomes fixed in the new coordinate, and we introduce the notation
\begin{align}
D_t = \frac{\p}{\p t}|_{y=\text{constant}} = \frac{\p}{\p t}|_{x=\text{constant}} + v^k\frac{\p}{\p x^k}. \label{lag}
\end{align}
to be the material derivative and
\begin{align*}
\p_i = \frac{\p}{\p x^i} = \frac{\p y^a}{\p x^i}\frac{\p}{\p y^a}.
\end{align*}
Due to (\ref{lag}), we shall also call $D_t$ as the time derivative as well by slightly abuse of terminology.\\
\indent Sometimes it is convenient to work in the Eulerian coordinate $(t,x)$, and sometimes it is easier to work in the Lagrangian coordinate $(t,y)$. In the Lagrangian coordinate the partial derivative $\p_t=D_t$ has more direct significance than it in the Eulerian frame. However, this is not true for spatial derivatives $\p_i$. The notion of space derivative that plays a more significant role in the Lagrangian coordinate is that the covariant differentiation with respect to the metric $g_{ab}(t,y)=\delta_{ij}\frac{\p x^i}{\p y^a}\frac{\p x^j}{\p y^b}$. We shall not involve covariant derivatives in our energy; instead, we use the regular Eulerian spatial derivatives. We will work mostly in the Lagrangian coordinate in this paper. However, our statements are coordinate independent.\\
\indent The Euclidean metric $\delta_{ij}$ in $\DD_{t}$ induces a metric
\begin{align}
g_{ab}(t,y)=\delta_{ij}\frac{\p x^i}{\p y^a}\frac{\p x^j}{\p y^b}, \label{g}
\end{align}
in $\Omega$ for each fixed $t$. We will denote covariant differentiation in the $y_{a}$-coordinate by $\nab_a$, $a=1,\cdots,n$, and the differentiation in the $x_i$-coordinate by $\p_i$, $i=1,\cdots,n$. Here, we use the convention that differentiation with respect to Eulerian coordinates is denoted by letters $i,j,k,l$ and with respect to Lagrangian coordinate is denoted by $a,b,c,d$.

\indent The regularity of the boundary is measured by the regularity of the normal, let $N^a$ to be the unit normal to $\p\Omega$,
$$
g_{ab}N^aN^b=1,
$$
and let $N_a=g_{ab}N^b$ denote the unit co-normal, $g^{ab}N_aN_b=1$. The induced metric $\gamma$ on the tangent space to the boundary $T(\p\Omega)$ extended to be $0$ on the orthogonal complement in $T(\Omega)$ is given by
$$
\gamma_{ab}=g_{ab}-N_aN_b,\q \gamma^{ab}=g^{ac}g^{bd}\gamma_{cd}=g^{ab}-N^aN^b.
$$
The orthogonal projection of an $(0,r)$ tensor $S$ to the boundary is given by
$$
(\Pi S)_{a_1,\cdots,a_r}=\gamma_{a_1}^{b_1}\cdots\gamma_{a_r}^{b_r}S_{b_1,\cdots,b_r},
$$
where $\gamma_{a}^{b}=g^{bc}\gamma_{ac}=\delta_{a}^{b}-N_aN^b$. In particular, the covariant differentiation on the boundary $\overline{\nab}$ is given by
$$
\overline{\nab}S=\Pi \nab S.
$$
We note that $\overline{\nab}$ is invariantly defined since the projection and $\nab$ are. The second fundamental form of the boundary $\theta$ is given by $\theta_{ab}=(\cnab N)_{ab}$, and the mean curvature of the boundary $\sigma=tr\theta=g^{ab}\theta_{ab}$.

It is now important to compute time derivative of the metric $D_tg$, the normal $D_tN$, as well as the time derivative of corresponding measures.

\lem \label{Lemma 1.1}
Let $x=f_t(y)=x(t,y)$ be the change of variable given by
\begin{equation}\label{lagrangiancoordinates}
\frac{dx}{dt}=v(t,x(t,y)), \q x(0,y)=f_0(y),\q y\in\Omega,
\end{equation}
and
\begin{equation*}
g_{ab}(t,y)=\delta_{ij}\frac{\p x^i}{\p y^a}\frac{\p x^j}{\p y^b},
\end{equation*}
to be the induced metric.
In addition, we let $\gamma_{ab} = g_{ab}-N_a N_b$, where $N_a= g_{ab}N^b$ is the co-normal to $\p\Omega$. Set
\begin{align}
v_a(t,y) = v_i(t,x)\frac{\p x^i}{\p y^a}, \q u^a = g^{ab}u_b,\\
d\mu_g, \q \text{volume element with respect to the metric}\,\, g ,\\
d\mu_\gamma,\q\text{surface element with respect to the metric}\,\, \gamma.\label{dmu}
\end{align}
Then
\begin{align}
D_t g_{ab} = \nab_a v_b+\nab_b v_a,\label{Dtg}\\
D_t g^{ab} = -g^{ac}g^{bd}D_tg_{cd},\label{Dtg_inverse}\\
D_t N_a = -\frac{1}{2}N_a(D_tg^{cd})N_cN_d,\label{DtN}\\
D_t d\mu_g = \di v\,d\mu_g\label{dg},\\
D_t d\mu_\gamma =(\sigma v \cdot N)\,d\mu_\gamma. \label{T2'}
\end{align}
\begin{proof} We have, since the commutator $[D_t,\frac{\p}{\p y}]=0$, and $D_t x(t,y) = v(t,y)$,
$$
D_t\frac{\p x^i}{\p y^a} = \frac{\p v_i}{\p y^a} = \frac{\p x^k}{\p y^a}\frac{\p v_i}{\p x^k},
$$
and so
\begin{align*}
D_t g_{ab} = \sum_{i}D_t(\frac{\p x^i}{\p y^a}\frac{\p x^i}{\p y^b})= \frac{\p x^k}{\p y^a}\frac{\p v_i}{\p x^k}\frac{\p x^i}{\p y^b}+\frac{\p x^i}{\p y^a}\frac{\p x^k}{\p y^b}\frac{\p v_i}{\p x^k} = \nab_a v_b+\nab_bv_a.
\end{align*}
(\ref{Dtg_inverse}) follows from $0=D_t(g^{ab}g_{bc})=D_t(g^{ab})g_{bc}+g^{ab}D_tg_{bc}$.
(\ref{dg}) follows since in local coordinate we have $d\mu_{g} = \sqrt{\detg}\,dy$ and $D_t\detg = \detg g^{ab}D_t g_{ab} = 2{\detg}\, \di v$.\\
To prove (\ref{DtN}), we choose the local foliation $u$ so that $\p\Omega=\{y:u(y)=0\}$ and $u<0$ in $\Omega$, then
$$
N_a=\frac{\p_a u}{\sqrt{{g^{cd}}\p_cu\p_du}},
$$
and (\ref{DtN}) follows from a direct computation. Now since
$$
d\mu_\gamma = \frac{\sqrt{\det g}}{\sqrt{\sum N_n^{2}}}\,dS(y),
$$
where $dS(y)$ is the Euclidean surface measure. By (\ref{DtN}) we have
$$
D_t d\mu_\gamma = \di v+\frac{1}{2}(D_tg^{cd})N_cN_d.
$$
But since $\di v=g^{ab}D_tg_{ab}/2$ and then (\ref{Dtg}) and (\ref{Dtg_inverse}) imply
$$
D_t d\mu_\gamma = \frac{1}{2}g^{ab}D_tg_{ab}-\frac{1}{2}(D_tg_{ab})N^aN^b=\gamma^{ab}\nab_av_b.
$$
But since $\gamma^{ab}\nab_av_b=\gamma^{ab}\cnab_a(N_bv\cdot N)+\gamma^{ab}\cnab_a\overline{v}_b$, and $\gamma^{ab}\cnab_a\overline{v}_b=\di v|_{\p\Omega}=0$, (\ref{T2'}) follows.
\end{proof}

\section{Estimates on a bounded domain with a moving boundary}\label{section 3}
Most of the results in this section will be stated in a coordinate-independent fashion. Throughout this section, $\nab$ will refer to covariant derivative with respect to the metric $g_{ij}$ in $\Omega$, and $\cnab$ will refer to covariant differentiation on $\p\Omega$ with respect to the induced metric $\gamma_{ij}=g_{ij}-N_iN_j$. Hence, in this section (and only), $\Omega$ will be used to denote a general domain with smooth boundary. In addition, we shall assume the normal $N$ to $\p\Omega$ is extended to a vector field in the interior of $\Omega$ satisfying $g_{ij}N^iN^j\leq 1$ by the same way introduced in Lemma \ref{trace 1}.
\subsection{Elliptic estimates}
\mydef\label{tensor} Let $u:\Omega\subset\RR^n\to\RR^n$ be a smooth vector field, and $\beta_k=\beta_{Ik}=\nab_{I}^{r}u_k$ be the $(0,r)$-tensor defined based on $u_k$, where $\nab_{I}^{r}=\nab_{i_1}\cdots\nab_{i_r}$ and $I=(i_1,\cdots,i_r)$ is the set of indices. Let $\di \beta_k=\nab_i\beta^{i}=\nab^r\di u$ and $\curl \beta=\nab_i\beta_j-\nab_j\beta_i=\nab^r\curl u_{ij}$.
\mydef (norms) If $|I|=|J|=r$, let $g^{IJ}=g^{i_1j_1}\cdots g^{i_r j_r}$ and $\gamma^{IJ}=\gamma^{i_1j_1}\cdots\gamma^{i_rj_r}$. If $\alpha$, $\beta$ are $(0,r)$ tensors, let $\langle\alpha,\beta\rangle=g^{IJ}\alpha_I \beta_J$ and $|\alpha|=\langle\alpha,\alpha\rangle$. If $(\Pi\beta)_{I}=\gamma_{I}^{J}\beta_J$ is the projection, then $\langle \Pi \alpha,\Pi\beta\rangle=\gamma^{IJ}\alpha_I\beta_J$.Let
\begin{align*}
||\beta||_{L^2(\Omega)}=(\int_{\Omega}|\beta|^2\,d\mu_g)^{\frac{1}{2}},\\
||\beta||_{L^2(\p\Omega)}=(\int_{\p\Omega}|\beta|^2\,d\mu_\gamma)^{\frac{1}{2}},\\
||\Pi\beta||_{L^2(\p\Omega)}=(\int_{\p\Omega}|\Pi\beta|^2\,d\mu_\gamma)^{\frac{1}{2}}.
\end{align*}
\indent We now state the following Hodge-type decomposition theorem, which serves as a main ingredient for proving the elliptic estimates.
\thm \label{hodge}(Hodge-decomposition) Let $\beta$ be defined in Definition \ref{tensor}. If $|\theta|+|\frac{1}{l_0}|\leq K$, where $\theta$ is the second fundamental form and $l_0$ is the injective radius defined in (\ref{inj rad}), then
\begin{align}
|\nab\beta|^2\lesssim g^{ij}\gamma^{kl}\gamma^{IJ}\nab_k\beta_{Ii}\nab_l\beta_{Jj}+|\di \beta|^2+|\curl\beta|^2\\
\int_{\Omega}|\nab\beta|^2\,d\mu_g\lesssim\int_{\Omega}(N^iN^jg^{kl}\gamma^{IJ}\nab_k\beta_{Ii}\nab_l\beta_{Jj}+|\di \beta|^2+|\curl\beta|^2+K^2|\beta|^2)\,d\mu_g.
\end{align}
\begin{proof}
See \cite{CL}
\end{proof}

\lem \label{poincare}(Poincar\'{e} type inequalities)
Let $q:\Omega\subset\RR^n\to\RR$ be a smooth and $q|_{\p\Omega}=0$, then
\begin{align}
||q||_{\lli}\lesssim (\vol\Omega)^{\frac{1}{n}}||\nab q||_{\lli},\label{FB}\\
||\nab q||\lesssim (\vol\Omega)^{\frac{1}{n}}||\lap q||_{\lli}.\label{Poincare}
\end{align}
\begin{proof} The first inequality is Faber-Krahns theorem, whose proof can be found in \cite{YS}. The second inequality follows from the first and integration by parts.
\end{proof}

\prop \label{ellpitic estimate}(Elliptic estimates) Let $q:\Omega\to\RR$ be a smooth function. Suppose that $|\theta|+|\frac{1}{l_0}|\leq K$, then we have, for any $r\geq2$ and $\delta>0$,
\begin{align}
||\nab^r q||_{\lli}+||\nab^r q||_{\llb}\lesssim_{K,\vol\Omega}\sum_{s\leq r}||\Pi\nab^s q||_{\llb}+\sum_{s\leq r-1}||\nab^s\lap q||_{\lli},\label{ell est I}\\
||\nab^r q||_{\lli}+||\nab^{r-1} q||_{\llb} \lesssim_{K,\vol\Omega}\delta\sum_{s\leq r}||\Pi\nab^s q||_{\llb}+\delta^{-1}\sum_{s\leq r-2}||\nab^s\lap q||_{\lli}.\label{ell est II}
\end{align}
where we have applied the convention that $A\lesssim_{p,q}B$ means $A\leq C_{p,q}B$.
\begin{proof}
See \cite{CL} (Proposition 5.8).
\end{proof}

\subsection{Estimate for the projection of a tensor to the tangent space of the boundary}
The use of the projection of the tensor $\Pi\nab^sD_t^kh$ in the boundary part of energy (\ref{Er}) is essential to compensate the potential loss of regularity. A simple observation that will help us is that if $q=0$ on $\p\Omega$, then $\Pi\nab^2q$ contains only first-order derivative of $q$ and all components of the second fundamental form. To be more precise, we have
\begin{equation}
\Pi\nab^2q=\cnab^2q+\theta\nab_Nq,\label{21}
\end{equation}
where the tangential component $\cnab^2q=0$ on the boundary. Furthermore, in $L^2$ norms, (\ref{21}) yields,
\begin{equation}
||\Pi\nab^2q||_{\llb}\leq |\theta|_{L^{\infty}(\p\Omega)}||\nab_Nq||_{\llb}.\label{22}
\end{equation}
To prove (\ref{21}), we first recall the components of the projection operator $\gamma_{i}^{j}=\delta_{i}^{j}-N_iN^j$, hence
$$
\gamma_{j}^{k}\nab_i\gamma_{k}^{l}=-\gamma_{j}^{k}\nab_i(N_kN^l)=-\gamma_j^k\theta_{ik}N^l-\gamma_j^kN_k\theta_i^l=-\theta_{ij}N^l,
$$
and so
\begin{dmath*}
\cnab_i\cnab_jq=\gamma_{i}^{i'}\gamma_{j}^{j'}\nab_{i'}\gamma_{j'}^{j''}\nab_{j''}q
=\gamma_i^{i'}\gamma_j^{j'}\gamma_{j'}^{j''}\nab_{i'}\nab_{j''}q+\gamma_i^{i'}\gamma_{j}^{j'}(\nab_{i'}\gamma_{j'}^{j''})\nab_{j''}q
=\gamma_{i}^{i'}\gamma_{j}^{j'}\nab_{i'}\nab_{j'}q-\theta_{ij}\nab_Nq.
\end{dmath*}
In general, the higher order projection formula is of the form
$$
\Pi \nab^r q = (\cnab^{r-2}\theta)\nab_N q+O(\nab^{r-1}q)+O(\cnab^{r-3}\theta),
$$
which suggests the following generalization of (\ref{22}), its detailed proof can be found in \cite{CL}.
\prop \label{tensor estimate}(Tensor estimate)
Suppose that $|\theta|+|\frac{1}{l_0}|\leq K$, and for $q=0$ on $\p\Omega$, then for $m=0,1$
\begin{multline}
||\Pi\nab^{r}q||_{\llb}\lesssim_{K} ||(\cnab^{r-2}\theta)\nab_{N}q||_{\llb}+\sum_{l=1}^{r-1}||\nab^{r-l}q||_{\llb},\\
 +(||\theta||_{\linf}+\sum_{0\leq l\leq r-2-m}||\cnab^{l}\theta||_{\llb})(\sum_{0\leq l\leq r-2+m}||\nab^l q||_{\llb}), \label{tensor est}
\end{multline}
where the second line drops for $0\leq r\leq 4$.
\begin{proof}
See \cite{CL} (Proposition 5.9).
\end{proof}
\subsection{Estimate for the second fundamental form}
The estimate of the second fundamental form is a direct consequence of Proposition \ref{tensor estimate} with $q=h$ together with the Taylor sign condition, e.g., $|\nab_N h|\geq \epsilon>0$.
\prop \label{theta estimate}($\theta$ estimate) \footnotemark
Assume that $0\leq r\leq4$. Suppose that $|\theta|+|\frac{1}{l_0}|\leq K$, and the Taylor sign condition $|\nab_N h|\geq \epsilon>0$ holds, then
\begin{equation}
||\cnab^{r-2}\theta||_{\llb}\lesssim_{K, \frac{1}{\epsilon}}||\Pi \nab^rh||_{\llb}+\sum_{s=1}^{r-1}||\nab^{r-s}h||_{\llb}. \label{theta est}
\end{equation}
In fact, (\ref{theta est}) can both be applied to the cases when $r>4$ by modifying the lower order terms. We refer \cite{CL} for the details.
\footnotetext{The $\theta$ estimate suggests that the boundary regularity is in fact controlled by the boundary $L^2$ -norm of $h$, with a loss of $2$ derivatives.}

\section{Energy estimates for the wave equation}\label{section 4}
In this section we study the estimates for the enthalpy $h$. The commutator between $D_t$ and $\p$ is of the form
\begin{align}
[D_t,\p_i] = -(\p_iv^k)\p_k. \label{Dtdi comm}
\end{align}
If we take divergence on the first equation of (\ref{EE}), together with the fact that $\di v=-D_te(h)$ and (\ref{Dtdi comm}), we obtain
\begin{align}
D_t^2e(h)-\lap h=(\p_iv^j)(\p_jv^i),\q \text{in}\q [0,T]\times \Omega, \label{WW}
\end{align}
with initial and boundary conditions
\begin{align}
h|_{t=0} = h_0,\q D_th|_{t=0} = h_1,
\end{align}
and
\begin{align}
h|_{\p\Omega}=0.
\end{align}
Here, $\lap h = \delta^{ij}\p_i\p_jh= \frac{1}{\sqrt{|\det g|}}\p_a(\sqrt{|\det g|}g^{ab}\p_bh)$.
\subsection{Some commutators}
We are able to obtain a higher order version of (\ref{WW}) by commutating more time derivatives to it. But since our $D_t$ no longer commutes with the spatial derivatives, we need to compute the following commutators first:
\begin{enumerate}
\item $[\p_i, D_t^k]=\sum_{l=0}^{k-1}D_t^l[\p_i,D_t]D_t^{k-l-1}$
\item $[\lap,D_t]=\lap v^j\p_j+2\p^iv^j\p_i\p_j=-\p^j(D_te(h))\p_j+\p^j\curl_{kj}v+2\p^iv^j\p_i\p_j$, where $\p^i = \delta^{ik}\p_k$
\\
\\
The second equality is because $\lap v_j = \sum_k \p_k\p_kv_j=\p_j\di v+\sum_k\p_k\curl_{kj}v$
\item $[\lap,D_t^{r-1}] = \sum_{l=0}^{r-2}D_t^l[\lap,D_t]D_t^{r-l-2}$
\end{enumerate}
\indent Although $D_t$ and $\p$ are not commutative, (\ref{Dtdi comm}) implies that the commutator between $D_t$ and $\p$ is free from time derivative. In general, $[D_t^k,\p]$ is a product of mixed space-time derivative where each component depends on at most $k-1$ time derivatives. In fact for $k\leq 4$, each component of $[\p,D_t^k]$ and $[\lap, D_t^k]$ can be controlled by the a priori assumptions. This can be seen by the simplified version of the commutators, by expressing them in the format of $\text{main terms}+\text{lower order terms}$. To do it, we would like to introduce the following short-hand notations first.
\mydef (Symmetric dot product)  Let $[D_t,\p]=-(\p v)\symdot\p$, where the symmetric dot product $(\p v)\symdot\p$ is define component-wise by $((\p v)\symdot\p)_i=\p_i v^k \p_k$. In general, we have
\begin{equation}
[D_t,\p^r]=\sum_{s=0}^{r-1}\p^s[D_t,\p]\p^{r-s-1}=\sum_{s=0}^{r-1}-
\begin{pmatrix}
r\\
s+1
\end{pmatrix}
(\p^{1+s}v) \symdot \p^{r-s},  \label{Dtd^r comm}
\end{equation}
where
$$
((\p^{1+s}v) \symdot \p^{r-s})_{i_1,\cdots,i_r}=\frac{1}{r!}\sum_{\sigma\in S_r}(\p^{1+s}_{i_{\sigma_1}\cdots i_{\sigma_{1+s}}}v^k)(\p^s_{k,i_{\sigma_{s+2}\cdots i_{\sigma_{r}}}}),
$$
where $S_r$ is the $r$-symmetric group.

Now, the commutators $[\p, D_t^k], k\geq 2$ and $[\lap, D_t^{r-1}],r\geq 3$ can be rewritten as
\begin{align}
[\p, D_t^k] = \sum_{l_1+l_2=k-1}c_{l_1,l_2}(\p D_t^{l_1}v)\symdot(\p D_t^{l_2})+ \sum_{l_1+\cdots+ l_n= k-n+1, \, n\geq 3} d_{l_1,\cdots,l_n}(\p D_t^{l_1} v)\cdots (\p D_t^{l_{n-1}} v) (\p D_t^{l_n}),\label{dDt^k comm}
\end{align}
and
\begin{multline}
[\lap, D_t^{r-1}] = \sum_{l_1+l_2=r-2}c_{l_1,l_2} (\lap D_t^{l_1}v)\cdot (\p D_t^{l_2})+ \sum_{l_1+l_2=r-2}c_{l_1,l_2} (\p D_t^{l_1}v)\cdot (\p^2 D_t^{l_2})\\
+ \sum_{l_1+\cdots+ l_n= r-n, \, n\geq 3} d_{l_1,\cdots,l_n}(\p D_t^{l_3} v)\cdots (\p D_t^{l_{n}} v)\cdot (\lap D_t^{l_1} v)\cdot(\p D_t^{l_2})\\
+ \sum_{l_1+\cdots+ l_n= r-n, \, n\geq 3} e_{l_1,\cdots,l_n}(\p D_t^{l_3} v)\cdots (\p D_t^{l_{n}} v)\cdot (\p^2 D_t^{l_1} v)\cdot(\p D_t^{l_2})\\
+\sum_{l_1+\cdots+ l_n= r-n, \, n\geq 3} f_{l_1,\cdots,l_n}(\p D_t^{l_3} v)\cdots (\p D_t^{l_n} v)\cdot (\p D_t^{l_1} v)\cdot (\p^2 D_t^{l_2}), \label{lapDt comm}
\end{multline}
where the regular dot product is defined by the trace of the symmetric dot.
\subsection{The Energies $W_r(t)$}
By commutating $D_t^{r-1}$ on both sides of (\ref{WW}), we obtain the higher order wave equation
\begin{align}
e'(h)D_t^{r+1}h-\lap D_t^{r-1}h = f_r+g_r, \label{WWr}
\end{align}
where \begin{align}
f_r = D_t^{r-1}(\p v\cdot\p v)+[D_t^{r-1},\lap]h \label{f_r},
\end{align}
and $g_r$ is sum of terms of the form
\begin{align}
e^{m}(h)(D_t^{i_1}h)\cdots(D_t^{i_m}h), \q i_1+\cdots+i_m=r+1,\quad 1\leq i_1\leq \cdots \leq i_m\leq r .  \label{g_r}
\end{align}
Now, let us define the energy
\begin{align}
W_r(t) = \frac{1}{2}||\sqrt{e'(h)}D_t^rh||_{\lli}+\frac{1}{2}||\nab D_t^{r-1}h||_{\lli}. \label{Wr}
\end{align}
The following estimate holds for $W_r$:
\thm \label{energy est for wave eq} (Energy estimates for $W_r$)
Let $W_r$ be defined as (\ref{Wr}), then
\begin{align}
\frac{dW_r^2}{dt}\lesssim_M W_r^2+W_r(||f_r||_{\lli}+||g_r||_{\lli}).
\end{align}
\begin{proof}
We start by differentiating $||\sqrt{e'(h)}D_t^rh||_{\lli}^2$, and since $|e^{(r)}(h)|\leq c_0$, $|D_th|\leq M$, and $D_tJ=(\di v)J$, where $J=\sqrt{\det g}$, one has
\begin{dmath}
\frac{d}{dt}\frac{1}{2}\int_{\Omega}e'(h)D_t^rh\cdot D_t^rh \,J\,dy \lesssim \int_{\Omega}e'(h)D_t^{r+1}h\cdot D_t^rhJ\,dy+W_r^2 = \int_{\Omega}\lap D_t^{r-1}h\cdot D_t^rh\,J\,dy+\int_{\Omega}(f_r+g_r)D_t^r h\, J\,dy +W_r^2= -\int_{\Omega}\nab D_t^{r-1}h\cdot \nab D_t^rh\, J\,dy+\int_{\Omega}(f_r+g_r)D_t^r h \,J\,dy +W_r^2. \label{41}
\end{dmath}
Now by (\ref{Dtdi comm}), $\nab D_t^rh=D_t\nab D_t^{r-1}h+\nab v\symdot \nab D_t^{r-1}h$, and so that
\begin{dmath}
-\int_{\Omega}\nab D_t^{r-1}h\cdot \nab D_t^rh\, J\,dy = -\frac{1}{2}\frac{d}{dt}||\nab D_t^{r-1}h||_{\lli}^2 +\int_{\Omega}\nab D_t^{r-1}h(\nab v\symdot \nab D_t^{r-1}h)J\,dy, \label{42}
\end{dmath}
where the last term on the right hand side of (\ref{42}) is included in $W_r$ since $|\nab v|\leq M$.
\end{proof}

\subsection{Estimates for $||f_r||_{\lli}$}
Let us now analyse $||f_r||$ and $||g_r||$.
By adopting our notations used in (\ref{dDt^k comm})-(\ref{lapDt comm}), we are able to express $f_r$ as
\begin{multline}
f_r = \sum_{l_1+l_2=r-1}c_{l_1,l_2}(\nab D_t^{l_1}v)\cdot(\nab D_t^{l_2}v)+
 \sum_{l_1+l_2=r-2}d_{l_1,l_2}(\lap D_t^{l_1}v)\cdot (\nab D_t^{l_2}h)\\
  +\sum_{l_1+l_2=r-2}e_{l_1,l_2}(\nab D_t^{l_1}v)\cdot(\nab^2 D_t^{l_2}h) +\text{error terms}, \label{fr}
\end{multline}
where the "error terms" refer to the terms generated by the commutators, which are of the form
\begin{multline}
e_r=\sum_{l_1+\cdots+ l_n= r+1-n, \, n\geq 3} g_{l_1,\cdots,l_n}(\p D_t^{l_3} v)\cdots (\p D_t^{l_{n}} v)\cdot (\p D_t^{l_1} v)\cdot(\p D_t^{l_2}v)\\
+ \sum_{l_1+\cdots+ l_n= r-n, \, n\geq 3} e_{l_1,\cdots,l_n}(\p D_t^{l_3} v)\cdots (\p D_t^{l_{n}} v)\cdot (\p^2 D_t^{l_1} v)\cdot(\p D_t^{l_2} h)\\
+\sum_{l_1+\cdots+ l_n= r-n, \, n\geq 3} f_{l_1,\cdots,l_n}(\p D_t^{l_3} v)\cdots (\p D_t^{l_n} v)\cdot (\p D_t^{l_1} v)\cdot (\p^2 D_t^{l_2}h).
\end{multline}
We need to estimate $||f_r||_{\lli}$ for $r=2,3,4,5$. Since our estimates include mixed space-time derivatives, we would like to use the following more appealing notations.
\mydef(Mixed Sobolev norms)
let $u(t,\cdot):\RR^n \to \RR$ be a smooth function. We define
\begin{align*}
||u||_{r,\,0} = \sum_{s+k=r,\, k<r}||\nab^s D_t^ku||_{\lli},\\
||u||_r = ||u||_{r,0}+||\sqrt{e'(h)}D_t^rh||_{\lli}.
\end{align*}
It is worth mentioning that when $r\leq 5$, the error terms of $f_r$ generated by commutating $D_t$ and the spatial derivatives can be controlled linearly by mixed Sobolev norms of $v$ and $h$ of order at most $r-1$. On the other hand, in order to estimate the time derivative of $W_{r+1}$ in the next section, we have to make sure that the $r$-th order Sobolev norms in our estimates for $||f_r||_{\lli}$, $3\leq r\leq 5$ do not include $||\nab^r h||_{\lli}$ and $||v||_r$.
\subsubsection{When r=2}
The main terms involved in $f_2$ can be bounded by
\begin{align*}
 ||(\nab D_tv)(\nab v)||_{\lli}\leq ||\nab v||_{L^{\infty}}||\nab^2 h||_{\lli},\\
 ||(\lap v)(\nab h)||_{\lli}\leq  ||\nab h||_{L^{\infty}}||\nab^2 v||_{\lli},\\
 ||(\nab v)(\nab^2 h)||_{\lli}\leq |\nab v|_{L^{\infty}}||\nab^2 h||_{\lli}.
\end{align*}
Since the error terms in $f_2$ is of the form $\nab v\cdot\nab v\cdot\nab v$ , we get
$$||f_2||_{\lli}\lesssim_M ||\nab^2 v||_{\lli}+||\nab^2 h||_{\lli}+||\nab v||_{\lli}.$$

\subsubsection{When r=3}
The first and the third terms of $f_3$ can be bounded by
\begin{align*}
||(\nab D_t^2v)(\nab v)||_{\lli}+||(\nab D_tv)(\nab D_tv)||_{\lli}&\lesssim_M||\nab^2D_th||_{\lli}+||\nab^2 v||_{\lli}+||\nab^2 h||_{\lli},\\
||(\nab v)(\nab^2h)||_{\lli}+||(\nab D_tv)(\nab^2h)||_{\lli}&\lesssim_ M||\nab^2h||_{\lli},
\end{align*}
respectively.
To bound the second term, it is easy to see that by the wave equation (\ref{WWr}) and the fact that $|e^{(r)}(h)|\leq c_0$, we get
\begin{align*}
||(\lap D_tv)(\nab h)||_{\lli}=||(\nab\lap h)(\nab h)||_{\lli}\lesssim_M ||e'(h)\nab D_t^2h||_{\lli}+\sum_{j=1,2}||\nab^j v||_{\lli}+||h||_{2,0},
\end{align*}
and\footnotemark
$$
||(\lap v)(\nab D_th)||_{\lli}\leq |\nab D_th|_{L^{\infty}}||\nab^2v||_{\lli}\lesssim_M ||\nab^2 v||_{\lli}.
$$
\footnotetext{We must include $|\nab D_th|_{L^{\infty}}$ in our a priori assumptions. Since otherwise we would have to estimate $||(\lap v)(\nab D_th)||_{\lli}$ by interpolation, which contributes to $||\nab^3 v||_{\lli}$ and it is part of $E_3$. Hence we would lose one derivative when estimating the second order boundary $L^2$ norms.}

The higher order terms in $e_3$ are essentially bounded by the corresponding terms
in $f_r$, for $r\leq 3$, we just estimated times $|\nab v|_{L^\infty}$, apart from a term of the form
$\nab v\cdot\nab^2 v\cdot \nab h$ which can be estimated by $\|\nab^2 v\|_{L^2}$.
Hence,
$$||f_3||_{\lli}\lesssim_M ||\nab^2 D_th||_{\lli}+||e'(h)\nab D_t^2h||_{\lli}+||h||_{2,0}+\sum_{j=1,2}||\nab^j v||_{\lli}.$$

\subsubsection{When r=4}
The first term of $f_4$ can be bounded by
\begin{multline*}
\sum_{l_1+l_2=3}||c_{l_1,l_2}(\nab D_t^{l_1}v)(\nab D_t^{l_2}v)||_{\lli}\leq |\nab v|_{L^{\infty}}||\nab D_t^3v||_{\lli}+|\nab^2h|_{L^{\infty}}||\nab D_t^2v||_{\lei}\\
\lesssim_M||\nab^2D_t^2h||_{\lli}+||\nab^2D_th||_{\lli}+||\nab [D_t^2,\nab]h||_{\lli}+||\nab (\nab v\cdot \nab h)||_{\lli}\\
\lesssim_M ||\nab^2D_t^2h||_{\lli}+||\nab^2D_th||_{\lli}+\sum_{j=2,3}||h||_{j,0}+||\nab^2v||_{\lli}.
\end{multline*}
Whereas the third term can be bounded by
\begin{align*}
\sum_{l_1+l_2=2}||e_{l_1,l_2}(\nab D_t^{l_1}v)(\nab^2 D_t^{l_2}h)||_{\lli}\lesssim_M ||\nab D_t^2v||_{\lli}+||\nab^2 D_th||_{\lli}+||\nab^2D_t^2h||_{\lli}.
\end{align*}
To bound the second term, by interpolation (\ref{int interpolation}) we have
\begin{multline}
\sum_{l_1+l_2=2}||(\lap D_t^{l_1}v)\cdot (\nab D_t^{l_2}h)||_{\lli}\lesssim_{K,M} |\nab v|_{L^{\infty}}\sum_{j=1,2}||\nab^jD_t^2h||_{\lli}+|D_t^2h|_{L^{\infty}}\sum_{j=2,3}||\nab^jv||_{\lli}\\
+||\nab^3 D_th||_{\lli}+\sum_{j=2,3}(||\nab^j v||_{\lli}+||\nab^j h||_{\lli})\\
\lesssim_{K,M} ||\nab^2D_t^2h||_{\lli}+||\nab^3D_th||_{\lli}+
\sum_{j=2,3}(||h||_{j,0}+||\nab^j v||_{\lli}) \label{mod f_4}.
\end{multline}
Most of the terms in $e_4$ can be bounded by corresponding terms in $f_r$, for $r\leq 4$, and similar terms in $e_3$ times a priori assumptions, apart from terms of the form $\nab v\cdot\nab^2D_tv\cdot\nab h$, whose $L^2$ norm can be bounded by $||\nab^3h||_{\lli}$.

Therefore, we sum up and get
$$
||f_4||_{\lli}\lesssim_{K,M} ||\nab^3 D_th||_{\lli}+||\nab^2D_t^2h||_{\lli}+\sum_{j=2,3}||h||_{j,0}+\sum_{j=1,2,3}||\nab^j v||_{\lli}.
$$
\subsubsection{When $r=5$ and $n\leq 4$}
The first and the third terms of $f_5$  can be estimated by through similar method as above.
 \begin{multline*}
\sum_{l_1+l_2=4}||(\nab D_t^{l_1}v)(\nab D_t^{l_2}v)||_{\lli}+\sum_{l_1+l_2=3}||(\nab D_t^{l_1}v)(\nab^2 D_t^{l_2}h)||_{\lli}\\
\lesssim_{K,M} ||\nab^2D_t^3h||_{\lli}+\sum_{1\leq i\leq 4}||v||_{i,0}+\sum_{2\leq i\leq4}||h||_{i,0}.
\end{multline*}
We need the Sobolev lemma (\ref{interior sobolev}) to bound $\sum_{l_1+l_2=3}d_{l_1,l_2}||(\lap D_t^{l_1}v) (\nab D_t^{l_2}h)||_{\lli}$, whose terms are bounded by
\begin{equation}
 ||\lap v\cdot \nab D_t^3h||_{\lli}\lesssim_K (\sum_{j=2,3}||\nab^jv||_{\lli})(\sum_{j=1,2}||\nab^jD_t^3h||_{\lli}) \label{mod f_5},
\end{equation}
\begin{equation*}
||\lap D_t v\cdot\nab D_t^2h||_{\lli}\lesssim_{K} |\nab^2h|_{L^{\infty}}\sum_{j=1,2}||\nab^jD_t^2h||_{\lli}+|D_t^2h|_{L^{\infty}}\sum_{j=3,4}||\nab^jh||_{\lli},
\end{equation*}
\begin{equation*}
||\lap D_t^2v\cdot \nab D_th||_{\lli}\lesssim_M ||\nab^2D_t^2 v||_{\lli},
\end{equation*}
and
 \begin{multline*}||\lap D_t^3v\cdot\nab h||_{\lli}\lesssim_M ||\nab\lap D_t^2h||_{\lli}+||\lap[D_t^2,\nab]h||_{\lli}\\
\lesssim_{K,M} ||\nab^3D_t^2h||_{\lli}+\sum_{j=2,3}||v||_{i,0}+\sum_{j=3,4}||h||_{i,0},
\end{multline*}
respectively. Most of the terms in the error term $e_5$ are essentially bounded by corresponding terms in $f_r$, for $r\leq 5$, and similar terms in $e_3$ and $e_4$ times a priori assumptions, apart from the terms of the form $\nab v\cdot\nab^2D_t^2v\cdot\nab h$, which is estimated by $||\nab^2 D_t^2v||_{\lli}$.\\
\indent Hence\footnotemark
\begin{multline*}
||f_5||_{\lli}\lesssim_{K,M} ||\nab^3D_t^2h||_{\lli}+(\sum_{j=2,3}||\nab^j v||_{\lli})||\nab^2D_t^3h||_{\lli}\\
+\sum_{1\leq i\leq 4}||v||_{i,0}+(\sum_{j=2,3}||\nab^j v||_{\lli})||\nab D_t^3h||_{\lli}+\sum_{2\leq i\leq4}||h||_{i,0}.
\end{multline*}

\footnotetext{ It can be seen that our estimates started to loss linearity in the highest orders when $r\geq 5$.}
\subsection{Estimates for $||g_r||_{\lli}$ for $r\leq 4$.}
\indent Since $g_r$ is a sum of $e^{m}(h)(D_t^{i_1}h)\cdots(D_t^{i_m}h), \q i_1+\cdots+i_m=r+1$ and $1\leq i_1\leq\cdots\leq i_m\leq r$, one cannot apply interpolation inequalities on estimating $g_r$. But since terms of the form $||(D_th)^lD_t^kh||_{\lli}$ and $||(D_t^2h)^lD_t^kh||_{\lli}$ can be bounded by $||D_t^kh||_{\lli}$ times $|D_th|_{L^{\infty}}$ or $|D_t^2h|_{L^{\infty}}$. Further, since $|e^{(m)}(h)|\leq c_0\sqrt{e'(h)}$, we get for $r\leq 4$ that

\begin{align}
||g_r||_{\lli}\lesssim_{M,c_0}\sum_{i\leq r}||\sqrt{e'(h)} D_t^ih||_{\lli}. \label{g_r mod 1}
\end{align}
\subsection{Estimates for $||g_r||_{\lli}$ for $r=5$ and $n\leq 4$.}
The only difference for estimating $g_5$ is that it contains a quadratic term $e''(h)D_t^3h\cdot D_t^3h$, whose $L^2$ norm is bounded via Sobolev lemma (\ref{interior sobolev}). Hence,
$$
||e''(h)(D_t^3h)^2||\leq |D_t^3h|_{L^{\infty}}||\sqrt{e'(h)}D_t^3h||_{\lli},
$$
but
$$
|D_t^3h|_{L^{\infty}} \lesssim_{K,\vol\Omega} ||\nab^2 D_t^3h||_{\lli}+ ||\nab D_t^3h||_{\lli},
$$
where we have used the fact $||D_t^3h||_{\lli}\lesssim_{\vol\Omega} ||\nab D_t^3h||_{\lli}$ as a consequence of (\ref{FB}). Therefore, we conclude
\begin{align}
||g_5||_{\lli}\lesssim_{K,M,c_0,\vol\Omega} \sum_{j\leq 5} ||\sqrt{e'(h)}D_t^jh||_{\lli} +(||\nab^2 D_t^3h||_{\lli}+ ||\nab D_t^3h||_{\lli})||\sqrt{e'(h)}D_t^3h||_{\lli} .\label{g_r mod 2}
\end{align}
\section{Energy estimates for the Euler equations}\label{section 5}
We are now ready to prove our main theorem.
\prop \label{a priori estimate}Let $E_r$ be defined as (\ref{Er}), then there are continuous functions $C_r$  such that, for $t\in[0,T]$ and $r\leq 4$,
\begin{align}
|\frac{dE_r(t)}{dt}|\leq C_r(K,\frac{1}{\epsilon},M,c_0,\vol\DD_t,E_{r-1}^*)E_r^*(t), \label{main 1}
\end{align}
provided that the assumptions \eqref{econd} on $e(h)$ and the a priori bounds (\ref{geometry_bound})-(\ref{Dtth}) hold.

Our estimates will mostly be in the Lagrangian coordinates, but we shall compute the time derivative $\frac{d}{dt}E_r$ in Eulerian coordinate, since then we do not need to worry about the Christoffel symbols.
\subsection{Computing $\frac{d}{dt}E_r$}

We first compute $\frac{d}{dt}E_{s,k}(t)$, where $E_{s,k}$ is defined as \eqref{Esk}, when $s> 0$.
\begin{multline}
\frac{d}{dt}E_{s,k} =\frac{1}{2}\int_{\DD_t}\rho D_t(\delta^{ij}Q(\p^sD_t^kv_i,\p^sD_t^{k}v_j)\dx
+\frac{1}{2}\int_{\DD_t}\rho D_t(e'(h)Q(\p^sD_t^k h,\p^s D_t^kh))\dx\\
+\frac{1}{2}\int_{\p\DD_t}\rho D_t(Q(\p^sD_t^k h,\p^sD_t^k h)\nu)-Q(\p^sD_t^k h,\p^sD_t^k h)\nu(\sigma v\cdot N)+\rho Q(\p^sD_t^k h,\p^sD_t^k h)D_t\nu\,dS. \label{DtEr'}
\end{multline}
The estimates (\ref{extending nab n})-(\ref{Dtgamma}) together with a priori assumptions imply\footnote{We refer Section 5 of \cite{CL}  for the detailed proof} $$|D_tq^{ij}|\lesssim M, \q |\p q^{ij}|\lesssim M+K,\q  |\sigma v\cdot N|_{L^{\infty}(\p\Omega)}\lesssim K+M,$$
$$
|D_t\nu|_{L^{\infty}(\p\Omega)} = |D_t(-\nab_N h)^{-1}|_{L^{\infty}(\p\Omega)}\lesssim 1+\frac{1}{M},
$$
and
\begin{align}
D_t\gamma^{ij}=-2\gamma^{im}\gamma^{jn}(\frac{1}{2}D_tg_{mn}) \label{bdy q}.
\end{align}
Since $|D_tq^{ij}|\lesssim M$ in the interior and on the boundary $q^{ij}=\gamma^{ij}$,
and by \eqref{bdy q} $D_t \gamma$ is tangential, so  that (\ref{DtEr'}) can then be reduced to
\begin{multline}
\frac{d}{dt}E_{s,k} \leq \int_{\DD_t}\rho \delta^{ij}Q(D_t\p^sD_t^kv_i,\p^sD_t^{k}v_j)\dx+\int_{\DD_t}\rho e'(h)Q(D_t\p^sD_t^k h,\p^s D_t^kh)\dx\\
+\int_{\p\DD_t}\rho Q(D_t\p^sD_t^k h,\p^sD_t^k h)\nu\,dS+C(K,L,M)(E_{s,k}+||h||_{r,0}^2+||v||_{r,0}^2). \label{DtEr}
\end{multline}

Now, our commutators (\ref{Dtd^r comm}) and (\ref{dDt^k comm}) yield, since $D_tv_i = -\p_i h$,
\begin{align}
D_t\p^sD_t^k v_i = -\p^sD_t^k\p_i h +\sum_{0\leq m\leq s-1}c_{msk}(\p^{m+1}v)\symdot\p^{s-m}D_t^kv_i \label{comm1},\\
D_t\p^rh + (\p_jh)\p^rv^j = \p^rD_th+\sum_{0\leq m \leq r-2}d_{mr}(\p^{m+1}v)\symdot \p^{r-m}h\label{comm2},\\
D_t\p^sD_t^k h = \p^sD_t^{k+1}h+\sum_{0\leq m\leq s-1}d_{msk}(\p^{m+1}v)\symdot\p^{s-m}D_t^kh,\q \text{for}\q k\geq 1. \label{comm3}
\end{align}
We control the term $||(\p^{m+1}v)\symdot\p^{s-m}D_t^kv_i||_{\lei}$ in (\ref{comm1}) and $||(\p^{m+1}v)\symdot\p^{s-m}D_t^kh||_{\lei}$ in (\ref{comm3}) for $s+k=r$, $k<r$ and $r\leq 4$.
\begin{itemize}
\item  Since $r\leq 4$ and $k<r$ imply $k\leq 3$, the term $||(\p^{m+1}v)\symdot\p^{s-m}D_t^kh||_{\lei}$ can be bounded by
$$|D_t^kh|_{L^{\infty}}\sum_{j\leq s+1}||\p^j v||_{\lei}+|\p v|_{L^{\infty}}\sum_{j\leq s}||\p^jD_t^kh||_{\lli},
$$
when $k=1,2$. \\
For $k=3$, we have
$$
||(\p v)\symdot\p D_t^3h||_{\lei}\lesssim_M ||\p D_t^3h||_{\lei},
$$
and for $k=0$,
$$
||(\p^{m+1}v)\symdot \p^{r-m}h||_{\lei}\lesssim_K |\p v|_{L^{\infty}}\sum_{j\leq r}||\p^j h||_{\lei}+|\p h|_{L^{\infty}}\sum_{j\leq r} ||\p^j v||_{\lei}.
$$
\item The term $||(\p^{m+1}v)\symdot\p^{s-m}D_t^kv_i||_{\lei}$ with $k\geq 1$ can be re-expressed as
$$
||(\p^{m+1}v)\symdot \p^{s-m+1}D_t^{k-1}h||_{\lei}+||(\p^{m+1}v)\symdot \p^{s-m}[D_t^{k-1},\p]h||_{\lei},
$$
since $k-1\leq 2$, the second term is bounded by $\sum_{i\leq r-1}(||v||_{i,0}+||h||_{i,0})$ and the first can be bounded similarly as above.
\end{itemize}
The above anaylsis shows that the $L^2$ norm of the sum in (\ref{comm1})-(\ref{comm3})  contribute only to $||v||_{r,0}$ and $||h||_{r,0}$ with $r<4$. Hence,
\begin{multline}
\frac{d}{dt}E_{s,k} \leq -\int_{\DD_t}\rho (\delta^{ij}Q(\p^sD_t^kv_i,\p^sD_t^{k}\p_jh)\dx+\int_{\DD_t}\rho e'(h)Q(\p^sD_t^k h,\p^s D_t^{k+1} h)\dx\\
+\int_{\p\DD_t}\rho Q(\p^sD_t^k h,D_t\p^sD_t^k h)\nu\,dS
+ C(K,M)(||v||_{r,0}+||h||_{r,0})(\sum_{i\leq r}||v||_{i,0}+||h||_{i,0}). \label{DtEr''}
\end{multline}

Now (\ref{dDt^k comm}) yields,
\begin{multline}
||\p^sD_t^k\p h-\p^{s+1}D_t^{k}h||_{\lei} \lesssim  \sum_{l_1+l_2=k-1}||\p^s(\p D_t^{l_1}v\symdot\p D_t^{l_2}h)||_{\lei}\\
+ \sum_{l_1+\cdots+ l_n= k-n+1, \, n\geq 3}||\p^s(\p D_t^{l_1} v\cdots \p D_t^{l_{n-1}} v\symdot \p D_t^{l_n}h)||_{\lei}, \label{error terms}
\end{multline}
and the last two terms are bounded by $\sum_{i\leq r}(||h||_{i,0}+||v||_{i,0})$.
Therefore,

\begin{multline}
\frac{d}{dt}E_{s,k} \leq -\int_{\DD_t}\rho (\delta^{ij}Q(\p^sD_t^kv_i,\p_j\p^{s}D_t^{k}h)\dx+\int_{\DD_t}\rho e'(h)Q(\p^sD_t^k h,\p^{s} D_t^{k+1} h)\dx\\
+\int_{\p\DD_t}\rho Q(\p^sD_t^k h,D_t\p^sD_t^k h)\nu\,dS
+C(K,M)(||v||_{r,0}+||h||_{r,0})(\sum_{i\leq r}||v||_{i,0}+||h||_{i,0}). \label{DtEr final}
\end{multline}
\indent If we integrate by parts in the first term
\begin{multline}
\int_{\DD_t}\rho \delta^{ij}Q(\p^sD_t^k \p_iv_j,\p^sD_t^kh)\dx+\int_{\DD_t}\rho e'(h)Q(\p^sD_t^k h,\p^{s} D_t^{k+1} h)\dx\\
+\int_{\p\DD_t}\rho Q(\p^sD_t^k h,D_t\p^sD_t^k h-\nu^{-1}N_i\p^sD_t^kv^i)\nu\,dS
+C(K,M)(||v||_{r,0}+||h||_{r,0})(\sum_{i\leq r}||v||_{i,0}+||h||_{i,0}). \label{int by part}
\end{multline}
But since $\p^sD_t^{k+1}e(h)$ equals $e'(h)\p^s D_t^{k+1}h$ plus a sum of terms of the form
$$
e^{(m)}(h)(\p^{i_1}D_t^{j_1}h)\cdots(\p^{i_m}D_t^{j_m}h),
$$
where
$$
(i_1+j_1)+\cdots+(i_m+j_m)\leq r+1,\q 1\leq i_1+j_1\leq\cdots\leq i_m+j_m\leq r.
$$
Therefore,
\begin{multline}
\int_{\DD_t}\rho \delta^{ij}Q(\p^sD_t^k \p_iv_j,\p^sD_t^kh)\dx=\int_{\DD_t}\rho Q(\p^sD_t^k\di v,\p^sD_t^kh)\dx \\
= -\int_{\DD_t}\rho Q(\p^sD_t^{k+1}e(h),\p^sD_t^kh)\dx\\
\leq-\int_{\DD_t}\rho e'(h)Q(\p^sD_t^{k+1}h,\p^sD_t^kh)\dx
 +C(K,M)(\sum_{i\leq r}||h||_{r,0})^2,
\end{multline}
so the first integral in (\ref{int by part}) cancels with the second term. \\
We recall $\nu = -(\p_Nh)^{-1}$, so that $\nu^{-1}N_i=\p_ih$. Hence, the boundary term in (\ref{int by part}) becomes
\begin{align}
\sum_{k+s=r,s>0}\int_{\p\DD_t}\rho Q(\p^sD_t^k h,D_t\p^sD_t^k h+(\p_ih)(\p^sD_t^kv^i)\nu\,dS. \label{reduced boundary int}
\end{align}
Now, since (\ref{comm2}) and (\ref{comm3}), (\ref{reduced boundary int}) becomes sum of the boundary inner product of $\Pi\p^sD_t^kh$ and
\begin{align}
\Pi(D_t\p^rh + (\p_jh)\p^rv^j) = \Pi\p^rD_th+\sum_{0\leq m \leq r-2}d_{mr}\Pi((\p^{m+1}v)\symdot \p^{r-m}h)\label{bdy int s=r},\\
\Pi (D_t\p^sD_t^k h+(\p_ih)(\p^sD_t^kv^i)) = \Pi\p^sD_t^{k+1}h+\Pi(\p_ih)(\p^sD_t^kv^i)+\sum_{0\leq m\leq s-1}d_{mr}\Pi((\p^{m+1}v)\symdot\p^{s-m}D_t^kh), \label{bdy int s<r}
\end{align}
for $k=0$ and $k>0$, respectively.\\

In addition, when $s=0$,
\begin{dmath}
\frac{d}{dt}E_{0,r} \leq -\int_{\DD_t}\rho \delta^{ij}(D_t^{r}\p_i h)(D_t^r v_j)\dx+\int_{\DD_t}\rho e'(h)(D_t^{r+1}h)(D_t^rh)\dx + C(M)||\sqrt{e'(h)}D_t^rh||_{\lei}^2 \label{full time energy 1},
\end{dmath}
where we have used the fact that $e''(h)\leq c_0\sqrt{e'(h)}$. Furthermore, since
\begin{multline}
||D_t^r\p h-\p D_t^{r}h||_{\lei} \lesssim  \sum_{l_1+l_2=r-1}||\p D_t^{l_1}v\symdot\p D_t^{l_2}h||_{\lei}\\
+ \sum_{l_1+\cdots+ l_n= r-n+1, \, n\geq 3}||\p D_t^{l_1} v\cdots \p D_t^{l_{n-1}} v\symdot \p D_t^{l_n}h||_{\lei}
\end{multline}
(\ref{full time energy 1}) becomes, after integrating by parts on the first integral on the RHS of (\ref{full time energy 1}),
\begin{dmath}
\frac{d}{dt}E_{0,r} \leq \int_{\DD_t}\rho \delta^{ij}( D_t^{r} h)(D_t^r \di v)\dx+\int_{\DD_t}\rho e'(h)(D_t^{r+1}h)(D_t^rh)\dx + C(M)||\sqrt{e'(h)}D_t^rh||_{\lei}^2+C(M)\sum_{i\leq r}(||h||_{i,0}+||v||_{i,0})^2  \label{full time energy 2}.
\end{dmath}
But since
$$
D_t^{r}\di v = - D_t^{r+1} e(h) = -e'(h)D_t^{r+1}h - g_r,
$$
and because $||\sqrt{e'(h)}D_t^rh||_{\lei}$ is part of $||h||_r$, (\ref{full time energy 2}) becomes
\begin{align}
\frac{d}{dt}E_{0,r}\leq C(M)\sum_{i\leq r}(||h||_{i}+||v||_{i,0})^2.
\end{align}

Furthermore, let $K_r$ be defined as \eqref{K_r}, we have
\begin{align}
\frac{d}{dt}K_r = 2\int_{\DD_t}\rho |\p^{r-1}\curl v|\cdot|D_t\p^{r-1} \curl v|\dx.
\end{align}
But since the curl satisfies the equation
 $$D_t\curl_{ij}v=-(\p_iv^k)(\curl_{kj}v)+(\p_jv^k)(\curl_{ki}v),$$
then
\begin{dmath}
|D_t\p^{r-1}\curl v| \leq |\p^{r-1} D_t\curl v|+\sum_{0\leq m\leq r-2}e_{mr}(\p^{m+1}v)\symdot \p^{r-1-m}\,\curl v\\
 \lesssim \sum_{0\leq m\leq r-1}e_{mr}(\p^{m+1}v)\symdot \p^{r-1-m}\,\curl v.
\end{dmath}
The term $||(\p^{m+1}v)\symdot \p^{r-1-m}\,\curl v||_{\lei}$ can be bounded by
\begin{align}
|\p v|_{L^{\infty}}\sum_{j\leq r-1}||\p^{j}\curl v||_{\lei}+|\curl v|_{L^{\infty}}\sum_{j\leq r-1}||\p^{j+1}v||_{\lei}.
\end{align}

On the other hand,
\begin{dmath}
\frac{dW_{r+1}^2}{dt} = 2W_{r+1}\frac{dW_{r+1}}{dt} \lesssim W_{r+1}(W_{r+1}+||f_{r+1}||_{\lei}+||g_{r+1}||_{\lei}) \\
\lesssim E_r + \sqrt{E_r}(||f_{r+1}||_{\lei}+||g_{r+1}||_{\lei}).
\end{dmath}
This comes from the energy estimates for the wave equation, e.g., Theorem \ref{energy est for wave eq}.

\indent Therefore, we have proved:
\thm Let $E_r$ be defined as (\ref{Er}), for $r\leq 4$ we have
\begin{multline}
|\frac{dE_r}{dt}|\lesssim_{K,M}E_r+\sum_{k+s=r,k,s>0}\bigg(||\Pi\p^s D_t^kh||_{\leb} \Big(||\Pi\p^s D_t^{k+1}h||_{\leb}\\
+||\Pi(\p_ih)(\p^sD_t^kv^i)||_{\leb}+\sum_{0\leq m\leq s-1}||\Pi((\p^{m+1}v)\symdot\p^{s-m}D_t^kh)||_{\leb}\Big)\bigg)\\
+||\Pi\p^r h||_{\leb}\Big(||\Pi \p^r D_t h||_{\leb}+\sum_{0\leq m \leq r-2}||\Pi((\p^{m+1}v)\symdot \p^{r-m}h)||_{\leb}\Big)\\
+C(K,M)(\sum_{i\leq r}||v||_{i,0}+||h||_{i})^2+ \frac{dW_{r+1}^2}{dt}.
\end{multline}
\rmk (\ref{bdy int s=r}) is essential in our energy estimates since $\Pi\p^r v$ is cancelled by the commutator, since there is no way to control $\Pi\p^r v$ on the boundary due to the loss of regularity. However, we can control $\Pi\p^sD_t^kv$ in (\ref{bdy int s<r}) on the boundary for $k>0$ since it can be reduced to $||\Pi\p^{s+1}D_t^{k-1}h||_{\leb}$ modulo error terms, which can then be controlled by elliptic estimates.
\mydef(Mixed boundary Sobolev norm)
let $u(t,\cdot):\RR^n \to \RR$ be a smooth function. We define
\begin{align*}
\lee u\ree_r = \sum_{k+s=r}||\nab^sD_t^ku||_{\llb}.
\end{align*}
\indent Now, let us get back to Lagrangian coordinate. Based on the computation we have as well as (\ref{tensor interpolation}), controlling $\frac{d}{dt}$ requires to bound
$$
||v||_{r,0},||h||_{r},\sum_{j\leq r-1}||\nab^jv||_{\llb}, \lee h\ree_r,
$$
and
$$\sum_{k+s=r,s\geq2}||\Pi\nab^s D_t^{k+1}h||_{\llb}.
$$

\thm \label{the key lemma}
With the a priori assumptions (\ref{geometry_bound})-(\ref{Dtth}) hold, we have for $0\leq r\leq 4$,
\begin{align}
||v||_{r,0}+||h||_r \leq C(K,M,c_0,\vol \DD_t)\sqrt{E_r^*}\label{interior estimates}.
\end{align}
In addition to that,
\begin{align}
||D_th||_r+\lee h\ree _r \leq C_r(K,M,c_0,\frac{1}{\epsilon},\vol\Omega) \sqrt{E_r^*} \label{boundary estimate 1},\q 0\leq r\leq 3\\
||D_th||_4+\lee h\ree _4 \leq C_r(K,M,c_0,\frac{1}{\epsilon},\vol\Omega) (1+\sqrt{E_3^*})\sqrt{E_4^*} \label{boundary estimate 2 dim 3}.
\end{align}
We shall prove (\ref{interior estimates}) in Section \ref{section 5.2}, and (\ref{boundary estimate 1})-(\ref{boundary estimate 2 dim 3}) in Section \ref{section 5.3}.
\subsection{Interior estimates, bounds for $||v||_{r,0}$,$||h||_{r}$}\label{section 5.2}
Our strategy is to first apply Theorem \ref{hodge} to control $||v||_{r,0}$ in terms of the energies $E_r$ and $L^2$ norm of $h$, and then we will apply our elliptic estimate (\ref{ell est II}) to control $||h||_r$.  We shall only focus on $r=4$, since the other cases follow from a similar method. Now, since
\begin{align}
||v||_{4,0} \lesssim_M
||\nab^4 v||_{\lli}+||\nab^4h||_{\lli}+||\nab^3D_th||_{\lli}+||\nab^2 D_t^2h||_{\lli}+\sum_{1\leq i\leq 3}(||v||_{i,0}+||h||_{i,0})\label{v_4}.
\end{align}
So the terms of order $4$ except for $||\nab^4v||$ can be combined with $||h||_4$. Now, Theorem \ref{hodge} yields,
\begin{align}
||\nab^4 v||_{\lli}\lesssim \sqrt{E_4}+||\nab^3 \di v||_{\lli}.
\end{align}
We recall that $\di v = -e'(h)D_th$, hence
\begin{align}
||\nab^4 v||_{\lli}\lesssim \sqrt{E_4}+ \sum_{1\leq j\leq 4}||h||_{j,0} .
\end{align}

\indent To bound $||h||_4$, since (\ref{ell est II}) provides, for each $k,s$ that $k+s=r$,
\begin{align}
||\nab^sD_t^kh||_{\lli}\lesssim_{K,M,\vol\Omega}||\Pi\nab^sD_t^kh||_{\llb}+\sum_{0\leq j\leq s-2}||\nab^{j}\lap D_t^kh||_{\lli},
\end{align}
for $s\geq 2$ (and so $k\leq 2$). The term $||\Pi\nab^sD_t^kh||_{\llb}$ bounded by $(||\nab h||_{L^{\infty}(\p\Omega)}E_r)^{\frac{1}{2}}$, by the construction of $E_r$. Furthermore, by the wave equation (\ref{WWr})
\begin{multline}
\sum_{0\leq j\leq s-2, 2\leq s\leq 4, s+k=r}||\nab^{j}\lap D_t^kh||_{\lli}\lesssim\sum_{0\leq j\leq 2}||\nab^jD_t^2e(h)||_{\lli}\\
+\sum_{0\leq j\leq 1}||\nab^jD_t^3e(h)||_{\lli}+||D_t^4e(h)||_{\lli}
+\sum_{1\leq i\leq 3}(||v||_{i,0}+||h||_{i}).
\end{multline}
But since $|e^{(l)}(h)|\leq c_0e'(h)$, by the same way as we did to control $||g_r||_{\lli}$,
\begin{align}
\sum_{0\leq j\leq s-2, 2\leq s\leq 4, k+s=r}||\nab^{j}\lap D_t^kh||_{\lli}\lesssim_{K,M,c_0,\vol\Omega}||\nab^2D_t^2h||_{\lli}+W_4+\sum_{1\leq i\leq 3}(||v||_{i,0}+||h||_{i}),
\end{align}
and if we apply (\ref{ell est II}) again with $q=D_t^2h$ to $||\nab^2D_t^2h||_{\lli}$, we have
\begin{dmath*}
||\nab^2 D_t^2h||_{\lli} \lesssim_{K,M,c_0,\vol \Omega}||\Pi\nab^2 D_t^2h||_{\llb}+||\lap D_t^2h||_{\lli}\\
\lesssim_{K,M,c_0,\vol \Omega}\sqrt{E_4^*}+W_4+\sum_{1\leq i\leq 3}(||v||_{i,0}+||h||_{i}).
\end{dmath*}
In addition, under the inductive assumption that $\sum_{1\leq i\leq 3}(||v||_{i,0}+||h||_{i})$ is already bounded\footnote{We mention here that the lowest order terms, namely $||\nab h||$ and $||D_th||$, are bounded via Lemma \ref{poincare}.}  by $\sqrt{E_3^*}+W_3$, and together with (\ref{v_4}), we get
\begin{align}
||v||_{4,0}+||h||_4
\lesssim_{K,M, c_0,\vol \Omega} W_4^*+\sqrt{E_4^*}.
\end{align}
Now, if we apply similar analysis to $||v||_{r,0},||h||_r$ for $r=2,3$ and by induction, we get that for each $r\leq 4$,
\begin{align}
\sum_{1\leq i\leq r}(||v||_{i,0}+||h||_i)
\lesssim_{K,M,c_0,\vol \Omega}\sqrt{E_r^*}+W_r^* .
\end{align}
Now, since $W_r^{*}$ is part of the energy $\sqrt{E_{r-1}^*}$, we have proved (\ref{interior estimates}).

\subsection{Boundary estimates, bounds for $\sum_{j\leq r-1}||\nab^j v||_{\llb}, \lee h\ree_r$ and $||\cnab^{r-2}\theta||_{\llb}$}\label{section 5.3}
The control of $\sum_{j\leq r-1}||\nab^j v||_{\llb}$ follows directly form the estimate of $\sum_{j\leq r}||\nab^j v||_{\lli}$ by trace Theorem (Theorem \ref{trace theorem}) . On the other hand,  we shall not estimate $\lee h\ree_r$ alone; instead, we estimate\footnotemark $||D_th||_r+\lee h\ree_r$ by (\ref{ell est II}). This has to be done since we need to estimate $||f_{r+1}||_{\lei}$ and $||g_{r+1}||_{\lei}$ by $E_r$. We will first do the cases when $r=2,3$ in order to get a general idea.

\footnotetext{The reason that we use the norm $||D_th||_r$ instead of $||h||_{r+1}$ is because the latter involves $||\nab^{r+1} h||$ which, after applying the elliptic and tensor estimates, gives $||(\cnab^{r-1}\theta)\nab_N h||_{\llb}$ but $||\cnab^{r-1}\theta||_{\llb}$ can only be controlled by $E_{r+1}$.}
\subsubsection{When $r=2$}
We estimate the mixed boundary $L^2$ norm $\lee h\ree_2$ by (\ref{ell est I})
\begin{multline}
\lee h\ree_2 \lesssim_{K,M,\vol\Omega}||\Pi\nab^2 h||_{\llb}+\sum_{j\leq 1}||\nab^j\lap h||_{\lli}+||\lap D_th||_{\lli}\\
\lesssim_{K,M,c_0,\vol\Omega}\sqrt{E_2}+\sum_{j\leq1}||e'(h)\nab^jD_t^2h||_{\lli}+||e'(h)D_t^3h||_{\lli}+\sum_{i\leq 2}(||v||_{i,0}+||h||_{i,0})\\
\lesssim_{K,M,c_0,\vol\Omega} \sqrt{E_2}+W_3+\sum_{i\leq 2}(||v||_{i,0}+||h||_{i,0}),
\end{multline}
and by (\ref{ell est II}) we get, for each $\delta>0$ that
\begin{align}
||D_th||_2\lesssim_{K,M,c_0,\vol\Omega} \delta||\Pi\nab^2D_th||_{\llb}+\delta^{-1}||\lap D_th||_{\lli}+W_3.
\end{align}
Now if we combine the interior and boundary estimates, we have for $0<\delta< 1$ that
\begin{multline}
||D_th||_2+\lee h\ree_2\lesssim_{K,M,c_0,\vol \Omega}  \sqrt{E_2}+W_3+\delta||\Pi\nab^2D_th||_{\llb}+\delta^{-1}||\lap D_th||_{\lli}+\sum_{i\leq 2}(||v||_{i,0}+||h||_{i,0})\\
\lesssim_{K,M,c_0,\vol \Omega} \sqrt{E_2}+\delta||\Pi\nab^2D_th||_{\llb}+\delta^{-1}\Big(W_3+\sum_{i\leq 2}(||v||_{i,0}+||h||_{i,0})\Big)\label{D_th h 2}.
\end{multline}
Further, (\ref{interior estimates}) would imply
$$
W_3+\sum_{i\leq 2}(||v||_{i,0}+||h||_{i,0})\lesssim_{K,M,c_0,\vol\Omega}\sqrt{E_2^*}.
$$
Since by (\ref{tensor est}) we have
\begin{align}
\delta||\Pi\nab^2D_th||_{\llb}\lesssim_{K} \delta(|\nab_N D_th|_{\linf}||\theta||_{\llb}+\sum_{j\leq 1}||\nab^jD_th||_{\llb}).
\end{align}
Now if we take $\delta=\delta(K,M,\vol\Omega)$ to be sufficiently small, the last term on the RHS can be combined with $\lee h\ree_2$ on the left (since $D_th=0$ on $\p\Omega$). Since $||\theta||_{\llb}\leq \epsilon^{-1}||\Pi\nab^2h||_{\llb}$, and so the first term is part of $\sqrt{E_2}$. Therefore,
$$
||D_th||_2+\lee h\ree_2\lesssim_{K,M,c_0,\frac{1}{\epsilon},\vol \Omega}\sqrt{E_2^*}+W_3 \lesssim \sqrt{E_2^*},
$$
since $W_3$ is part of $\sqrt{E_2}$.

\subsubsection{When $r=3$}
By (\ref{ell est I}), we get
\begin{multline}
\lee h\ree_3 \lesssim_{K,M,\vol\Omega}\sum_{k+s=3,s>0}||\Pi\nab^s D_t^k h||_{\llb}+\sum_{j\leq 2}||\nab^j\lap h||_{\lli}+\sum_{j\leq 1}||\nab^j\lap D_th||_{\lli}+||\lap D_t^2h||_{\lli}\\
\lesssim_{K,M,c_0,\vol\Omega}\sqrt{E_3}+\sum_{j\leq2}||e'(h)\nab^j D_t^2h||_{\lli}+\sum_{j\leq 1}||e'(h)\nab^jD_t^3h||_{\lli}+||e'(h)D_t^4h||_{\lli}+\sum_{i\leq 3}(||v||_{i,0}+||h||_{i,0}),
\end{multline}
together with (\ref{interior estimates}) we have
\begin{align}
\lee h\ree_3 \lesssim_{K,M,\vol\Omega}\sqrt{E_3^*}+W_4^* +||\nab^2 D_t^2h||_{\lli},
\end{align}
where the last term is part of $||D_th||_3$.

On the other hand, by (\ref{ell est II}) with $0<\delta<1$ we get
\begin{multline}
||D_th||_3\lesssim_{K,M,\vol\Omega}\delta(||\Pi\nab^3D_th||_{\llb}+||\Pi\nab^2D_t^2h||_{\llb})\\
+\delta^{-1}(\sum_{j\leq 1}||\nab^j\lap D_th||_{\lli}+||\lap D_t^2h||_{\lli})+W_4\\
\lesssim_{K,M,c_0,\vol\Omega}\delta(||\Pi\nab^3D_th||_{\llb}+||\Pi\nab^2D_t^2h||_{\llb})
+\delta^{-1}(W_4+\sum_{i\leq 3}(||v||_{i,0}+||h||_{i,0})).
\end{multline}
Now (\ref{tensor est}) implies
\begin{dmath}
\delta(||\Pi\nab^3D_th||_{\llb}+||\Pi\nab^2D_t^2h||_{\llb})\lesssim_K \\
\delta(||\cnab \theta||_{\llb}|\nab_ND_th|_{\linf}+|\theta|_{\linf}||\nab_ND_t^2h||_{\llb}\\
+\sum_{0\leq j\leq 2}||\nab^jD_th||_{\llb}+\sum_{0\leq j\leq1}||\nab^j D_t^2h||_{\llb}).
\end{dmath}
Now let $\delta$ to be sufficiently small, and so the last three terms on the second line can be absorbed into $\sum_{ j\leq 3}\lee h\ree_j$. In addition, since we have just proved that
$$
||\nab^2h||_{\llb}+||\nab D_th||_{\llb}\lesssim_{K,M,c_0, \vol \Omega}\sqrt{E_2^*},
$$
and so by (\ref{theta est}) we have
\begin{align}
||\cnab \theta||_{\llb}\lesssim_{K,\frac{1}{\epsilon}}\sqrt{E_3}+\sum_{1\leq j\leq 2}||\nab^jh||_{\llb}\lesssim_{K,M,c_0,\frac{1}{\epsilon}, \vol \Omega} \sqrt{E_3^*}.
\end{align}
Therefore, if we combine the estimates for $\lee h\ree_3$, $||D_th||_3$ and $||\cnab \theta||_{\llb}$, as well as the lower order $L^2$ norms, we get by (\ref{interior estimates}) that
\begin{align}
\sum_{1\leq i\leq 3}(||v||_{i,0}+||h||_i+\lee h\ree_i)+||D_th||_3\lesssim_{K,M,c_0,\frac{1}{\epsilon},\vol (\Omega)}W_{4}^*+\sqrt{E_3^*}.
\end{align}
Therefore, since $W_4^*$ is part of $\sqrt{E_3^*}$, we conclude
$$
||D_th||_3+\lee h\ree_3\lesssim_{K,M,c_0,\frac{1}{\epsilon},\vol (\Omega)}\sqrt{E_3^*}.
$$
\subsubsection{When $r=4$}
The estimates for $\lee h\ree_4$ and $||D_th||_4$ follows from the same analysis that we applied for the previous cases.
\begin{multline}
\lee h\ree_4\lesssim_{K,M,\vol \Omega}\sum_{k+s=4,s>0}||\Pi\nab^sD_t^kh||_{\llb}+\sum_{j\leq 3}||\nab^j\lap h||_{\lli}\\
+\sum_{j\leq 2}||\nab^j\lap D_th||_{\lli}+\sum_{j\leq 1}||\nab^j\lap D_t^2h||_{\lli}
+||\lap D_t^3h||_{\lli}\\
\lesssim_{K,M,c_0,\vol\Omega}\sqrt{E_4}+W_5^*+||\nab^2 D_t^3 h||_{\lli}+\sum_{j\leq 4}(||v||_{i,0}+||h||_{i,0}),
\end{multline}
and for $0<\delta<1$,
\begin{multline}
||D_th||_4\lesssim_{K,M,\vol\Omega}\delta(||\Pi\nab^4D_th||_{\lli}+||\Pi\nab^3D_t^2h||_{\lli}+||\Pi\nab^2D_t^3h||_{\lli})\\
+\delta^{-1}(\sum_{j\leq 2}||\nab^j\lap D_th||_{\lli}+\sum_{j\leq 1}||\nab^j\lap D_t^2h||_{\lli}+||\lap D_t^3h||_{\lli})+W_5\\
\lesssim_{K,M,c_0,\vol\Omega}\delta(||\Pi\nab^4D_th||_{\lli}+||\Pi\nab^3D_t^2h||_{\lli}+||\Pi\nab^2D_t^3h||_{\lli})
+\delta^{-1}(\sum_{j\leq 4}(||v||_{i,0}+||h||_{i,0})+W_5).
\end{multline}
The $L^2$ norm of the projected tensors can be estimated by
\begin{dmath}
\delta(||\Pi \nab^4D_th||_{\llb}+||\Pi \nab^2 D_t^3h||_{\llb})\\
\lesssim_{K}\delta(|\nab_ND_th|_{\linf}||\cnab^2\theta||_{\llb}+|\theta|_{\linf}||\nab_N D_t^3h||_{\llb}\\
+ \sum_{j\leq 3}||\nab^jD_th||_{\llb}+\sum_{j\leq 1}||\nab^jD_t^3h||_{\llb}),
\end{dmath}
In addition,
$$
\|\cnab^2\theta\|_{\lli}\lesssim_{K,\frac{1}{\epsilon}}||\Pi\nab^4h||_{\llb}+\sum_{i=1}^3||\nab^{i}h||_{\llb}\lesssim_{K,M,c_0,\frac{1}{\epsilon},\vol\Omega}\sqrt{E_4^*},
$$
and so $||\Pi \nab^4D_th||_{\llb}$ and $||\Pi \nab^2 D_t^3h||_{\llb}$ can be treated similarly as we did in the previous cases.
On the other hand,
\begin{align}
\delta||\Pi\nab^3D_t^2h||_{\llb}\lesssim_K \delta(||(\nab_ND_t^2h)\cnab\theta||_{\llb}+\sum_{ j\leq2}||\nab^jD_t^2h||_{\llb}).
\end{align}
The first term $||(\nab_ND_t^2h)\cnab\theta||_{\llb}$ is bounded via Gagliardo-Nirenberg interpolation inequality (Theorem \ref{gag-ni thm}) if $\Omega\in \RR^3$ (e.g., $\p\Omega\in \RR^2$),
\begin{multline}
||(\nab_ND_t^2h)\cnab\theta||_{\llb}\leq ||\nab_N D_t^2h||_{L^4}||\cnab \theta||_{L^4}\lesssim_K \|\nab_N D_t^2h\|_{L^2}^{\frac{1}{2}}\|\cnab\theta\|_{L^2}^{\frac{1}{2}}||\nab D_t^2h||_{H^1(\p\Omega)}^{\frac{1}{2}}||\cnab\theta||_{H^1(\p\Omega)}^{\frac{1}{2}}\\
\lesssim_{K,M,\frac{1}{\epsilon}, \vol\Omega} \sqrt{E_3^*}(||\nab D_t^2h||_{H^1(\p\Omega)}+||\cnab\theta||_{H^1(\p\Omega)})
\lesssim_{K,M,c_0,\frac{1}{\epsilon}, \vol\Omega}\sqrt{E_3^*}\sqrt{E_4^*}+\sqrt{E_3^*}||\nab D_t^2h||_{H^1(\p\Omega)},
\end{multline}
where the last term $||\nab D_t^2h||_{H^1(\p\Omega)}$ is part of $\lee h\ree_4$.\\
If $\Omega\in \RR^2$, we have
\begin{multline}
||(\nab_ND_t^2h)\cnab\theta||_{\llb}\leq |\nab_N D_t^2h|_{\linf}||\cnab\theta||_{\llb}\lesssim_K(\sum_{j\leq 2}||\nab^2 D_t^2h||_{\llb})||\cnab\theta||_{\llb}\\
\lesssim_{K,M,c_0,\frac{1}{\epsilon}, \vol\Omega}\sqrt{E_3^*}\,\,(\sum_{j\leq 2}||\nab^jD_t^2h||_{\llb}) \label{dim 2 without GN}.
\end{multline}

Now, if we combine the estimates for $\lee h\ree_4$, $||D_th||_4$ and $||\cnab^2 \theta||_{\llb}$,  as well as the lower order $L^2$ norms, we get
\begin{multline}
\sum_{1\leq i\leq 4}(||v||_{i,0}+||h||_i+\lee h\ree_i)+||D_th||_4
 \lesssim_{K,M,c_0,\frac{1}{\epsilon},\vol\Omega}  \delta^{-1}\sqrt{E_4^*}+ \delta\sqrt{E_3^*}\sqrt{E_4^*}+ \delta(|\theta|_{\linf}||\nab_N D_t^3h||_{\llb}\\
+ \sum_{j\leq 3}||\nab^jD_th||_{\llb}+\sum_{ j\leq2}||\nab^jD_t^2h||_{\llb}+\sum_{j\leq 1}||\nab^jD_t^3h||_{\llb}+\sqrt{E_3^*}||\nab D_t^2h||_{H^1(\p\Omega)}).
\end{multline}
Therefore, with $\delta$ chosen to be of the form $$\frac{C(K,M,c_0,\frac{1}{\epsilon},\vol\Omega)}{2(1+\sqrt{E_3^*})},$$ where $C$ is a continuous function that is sufficiently small, the above inequality implies
\begin{align}
\sum_{1\leq i\leq 4}(||v||_{i,0}+||h||_i+\lee h\ree_i)+||D_th||_4 \lesssim_{K,M,c_0,\frac{1}{\epsilon},\vol\Omega} (1+\sqrt{E_3^*})\sqrt{E_4^*},
\end{align}
and so
\begin{align}
||D_th||_4+\lee h\ree_4 \lesssim_{K,M,c_0,\frac{1}{\epsilon},\vol\Omega} (1+\sqrt{E_3^*})\sqrt{E_4^*}.
\end{align}

\subsection{Bounds for $\sum_{k+s=r,s\geq 2}||\Pi\nab^s D_t^{k+1}h||_{\llb}$ }
By (\ref{tensor est}), since $D_th=0$, we have
\begin{align}
\sum_{k+s=r,s\geq 2}||\Pi\nab^s D_t^{k+1}h||_{\llb}\lesssim_K \sum_{k+s=r,s\geq 2}||\cnab^{s-2}\theta\,(\nab_N D_t^{k+1}h)||_{\llb}\\
+\sum_{k+s=r, s\geq 2}\sum_{j=1}^{s-1}||\nab^{s-j}D_t^{k+1}h||_{\llb},
\end{align}
where the first term is bounded similarly by the arguments we had in the previous section. The second term is part of $\lee h\ree_r$.
\subsection{Bounds for $\frac{dW_{r+1}^2}{dt}$}\label{section 5.5}
We recall that we have
\begin{align}
\frac{dW_{r+1}^2}{dt} \lesssim E_r + \sqrt{E_r}(||f_{r+1}||_{\lli}+||g_{r+1}||_{\lli}).
\end{align}
Therefore, it suffices to bound $||f_{r+1}||_{\lli}$ and $||g_{r+1}||_{\lli}$ by $\sqrt{E_r^*}$, for $r=2,3,4$. On the other hand, we have proved in Section 4 that
\begin{align}
||f_3||_{\lli}\lesssim_M ||\nab^2 D_th||_{\lli}+||e'(h)\nab D_t^2h||_{\lli}+||h||_{2,0}+\sum_{j=1,2}||\nab^j v||_{\lli},\\
||f_4||_{\lli}\lesssim_{K,M} ||\nab^3 D_th||_{\lli}+||\nab^2D_t^2h||_{\lli}+\sum_{j=2,3}||h||_{j,0}+\sum_{j=1,2,3}||\nab^j v||_{\lli},
\end{align}
\begin{dmath}
||f_5||_{\lli}\lesssim_{K,M} ||\nab^3D_t^2h||_{\lli}+(\sum_{j=2,3}||\nab^j v||_{\lli})||\nab^2D_t^3h||_{\lli}\\
+\sum_{1\leq i\leq 4}||v||_{i,0}+(\sum_{j=2,3}||\nab^j v||_{\lli})||\nab D_t^3h||_{\lli}+\sum_{2\leq i\leq4}||h||_{i,0},
\end{dmath}
and
\begin{align}
||g_3||_{\lli}&\lesssim_{M,c_0} W_3^*,\\
||g_4||_{\lli}&\lesssim_{M,c_0} W_4^*,\\
||g_5||_{\lli}&\lesssim_{K,M,c_0,\vol\Omega}W_5^* +(||\nab^2 D_t^3h||_{\lli}+ ||\nab D_t^3h||_{\lli})||\sqrt{e'(h)}D_t^3h||_{\lli} .
\end{align}
Now, since $W_{r+1}^*$ is part of $\sqrt{E_r^*}$, (\ref{interior estimates}) and (\ref{boundary estimate 1})-(\ref{boundary estimate 2 dim 3}) give \footnote{We remark here that the estimates for $||f_r||$ do not require $||\nab^r h||_{\lli}$ norm. }
\begin{align}
||f_3||_{\lli}+||g_3||_{\lli} \lesssim_{K,M,c_0,\frac{1}{\epsilon},\vol \Omega} \sqrt{E_2^*}, \\
||f_4||_{\lli}+||g_4||_{\lli} \lesssim_{K,M,c_0,\frac{1}{\epsilon},\vol \Omega} \sqrt{E_3^*},
\end{align}
and finally
\begin{align}
||f_5||_{\lli}+||g_5||_{\lli} \lesssim_{K,M,c_0,\frac{1}{\epsilon},\vol \Omega} (1+\sqrt{E_3^*})^2\sqrt{E_4^*}. \q \text{when}\,\, \Omega\in\RR^3
\end{align}
\subsection{The energy estimates }
We are now ready to prove Proposition \ref{a priori estimate}. Since we have showed that our energies $E_r$ control the interior and boundary Sobolev norms of $v$ and $h$, the only thing left is to control the product of the projected tensors, i.e.,
\begin{align}
\sum_{s+k=r,s>0}\Big(\sum_{0\leq m\leq s-1}\Pi((\nab^{m+1}v)\symdot\nab^{s-m}D_t^kh)\Big),\q \text{for}\,\, k>0 \label{tensor k>0}\\
\sum_{0\leq m \leq r-2}\Pi((\nab^{m+1}v)\symdot \nab^{r-m}h),\q \text{for}\,\, k=0 \label{tensor k=0}
\end{align}
\begin{align}
\sum_{s+k=r,s>0}\Pi((\nab h)\symdot (\nab^sD_t^kv)).\q \text{for}\,\, k>0 \label{tensor special}
\end{align}

We cannot use interpolation (\ref{bdy interpolation}) here since it only applies to tangential derivative $\cnab$. Our strategy is to apply Sobolev lemma (Lemma \ref{boundary soboolev}) and Gagliardo-Nirenberg inequality to control terms that involving mixed derivatives\footnote{We want our estimates to be linear in the highest order. One can use Sobolev lemma only to control mixed Sobolev norms as well but the highest order energy would appear quadratically that way.}. Further, we use Theorem \ref{tensor interpolation} to control terms with full spatial derivatives.\\
\indent We bound (\ref{tensor k>0})-(\ref{tensor special}) when $r=4$ and $\Omega\in\RR^3$, since other cases follow from the same method and so we omit the details. By letting $\alpha=\nab^{s-1}v$ in (\ref{trace}) we get
$$
||\nab^{s-1}v||_{\llb}\lesssim_K \sum_{j\leq s}||\nab^j v||_{\lli}.
$$
Therefore, as we claimed in the beginning of Section 5.3, $\lee v\ree_{r-1}\lesssim_K \sum_{i\leq r}||v||_{i,0}$.

Now, each term of (\ref{tensor k>0}) is bounded as
\begin{itemize}
\item When $s=1,k=3$ (hence $m=0$)
$$||\Pi(\nab v\symdot\nab D_t^3h)||_{\llb}\lesssim_M ||\nab D_t^3h||_{\llb}\lesssim_{K,M,c_0,\frac{1}{\epsilon},\vol \Omega}(1+\sqrt{E_3^*})\sqrt{E_4^*}.$$
\item When $s=2, k=2$ (hence $m=0,1$)
\begin{dmath*}
||\Pi(\nab v\symdot \nab^2D_t^2h)||_{\llb}+||\Pi(\nab^2 v\symdot\nab D_t^2h)||_{\llb}
\\
\lesssim_{K,M}||\nab^2D_t^2h||_{\llb}+||\nab^2v||_{\llb}^{1/2}||\nab D_t^2h||_{\llb}^{1/2}||\nab^2v||_{H^1}^{1/2}||\nab D_t^2h||_{H^1}^{1/2}\\
\lesssim_{K,M,c_0,\frac{1}{\epsilon},\vol \Omega} (1+\sqrt{E_3^*})^2\sqrt{E_4^*}.
\end{dmath*}
\item When $s=3$,$k=1$ (hence $m=0,1,2$)
\begin{dmath*}
||\Pi(\nab v\symdot \nab^3 D_th)||_{\llb}+||\Pi(\nab^2 v\symdot \nab^2D_th)||_{\llb}+||\Pi(\nab^3 v\symdot \nab D_th)||_{\llb}\\
\lesssim_{K,M}||\nab^3 D_th||_{\llb}+||\nab^2v||_{\llb}^{1/2}||\nab^2 D_t h||_{\llb}^{1/2}||\nab^2v||_{H^1}^{1/2}||\nab^2 D_t h||_{H^1}^{1/2}+||\nab^3 v||_{\llb}\\
 \lesssim_{K,M,c_0,\frac{1}{\epsilon},\vol \Omega}(1+\sqrt{E_3^*})^2\sqrt{E_4^*}
\end{dmath*}
\end{itemize}
\rmk If $\Omega\in\RR^2$, the terms $||\Pi(\nab^2 v\symdot\nab D_t^2h)||_{\llb}$ and $||\Pi(\nab^2 v\symdot \nab^2D_th)||_{\llb}$ are bounded via Sobolev lemma. To be more specific,
\begin{dmath*}
||\Pi(\nab^2 v\symdot\nab D_t^2h)||_{\llb}+||\Pi(\nab^2 v\symdot \nab^2D_th)||_{\llb}\\
\lesssim_K (\sum_{j\leq2}||\nab^jD_t^2h||_{\llb})||\nab^2 v||_{\llb}+ (\sum_{j\leq3}||\nab^jD_th||_{\llb})||\nab^2 v||_{\llb}
\lesssim (1+ \sqrt{E_3^*})^2\sqrt{E_4^*}.
\end{dmath*}

The $L^2$ norm of (\ref{tensor k=0}) is bounded by the same arguments used in \cite{CL}. In addition, we need the following "product rule" of the projection, in order to control intermediate terms linearly in the highest order.
\thm  Let $S,R$ be two tensors, we have
\begin{align}
\Pi(S\symdot R)=\Pi(S)\symdot\Pi(R)+\Pi(S\symdot N)\tilde{\otimes} \Pi(N\symdot R), \label{proj prod rule}
\end{align}
where $\tilde{\otimes}$ is the symmetric tensor product which is defined similar to the symmetric dot product.
\begin{proof}
(\ref{proj prod rule}) is a direct consequence of the fact $g^{ab}=\gamma^{ab}+N^aN^b$.
\end{proof}
Now if we apply (\ref{proj prod rule}) to (\ref{tensor k=0}), we get for $0\leq m\leq r-2$
\begin{multline}
||\Pi(\nab^{1+m}v\symdot\nab^{r-m}h)||_{\llb}\leq \\||(\Pi\nab^{1+m}v)(\Pi\nab^{r-m}h)||_{\llb}
+||(\Pi N^j\nab^{1+m}v_j)(\Pi N^j\nab^{r-1-m}\nab_j h)||_{\llb}\\
\leq ||\Pi(\nab^{1+m}v)||_{L^{2(r-2)/m}(\p\Omega)}||\Pi\nab^{r-m}h||_{L^{2(r-2)/(r-2-m)}(\p\Omega)}\\ + ||\Pi(N^j\nab^{1+m}v_j)||_{L^{2(r-2)/m}(\p\Omega)}||\Pi N^j\nab^{r-m-1}\nab_jh||_{L^{2(r-2)/(r-2-m)}(\p\Omega)} \label{spatial tensor}.
\end{multline}
When $r=4$, if $m=0$ and $m=2$, this can be bounded via a priori assumptions time energies of the form $\sqrt{E_r^*}$, and the intermediate term $||\Pi(\nab^2v\symdot\nab^3h)||_{\llb}$ can be bounded linear in the highest orders by Theorem \ref{tensor interpolation}. The first term in the second line of (\ref{spatial tensor}) is bounded by letting $\alpha=\nab v$ and $\beta=\nab^2 h$ in (\ref{tensor int}). To be specific,
\begin{multline}
||(\Pi\nab^{2}v)(\Pi\nab^{3}h)||_{\llb}\lesssim_{K,M} (|\nab v|_{L^{\infty}}+\sum_{ j\leq 2}||\nab^jv||_{\llb})||\nab^4h||_{\llb}\\
+(|\nab^2 h|_{L^{\infty}}+\sum_{j\leq 3}||\nab^jh||_{\llb})||\nab^3 v||_{\llb}\\
+(||\theta||_{L^{\infty}}+||\cnab^2 \theta||_{\llb})(|\nab v|_{L^{\infty}}+\sum_{ j\leq 2}||\nab^jv||_{\llb})(|\nab^2 h|_{L^{\infty}}+\sum_{j\leq 3}||\nab^jh||_{\llb}).
\end{multline}
The term $||(\Pi N^j\nab^2 v_j)(\Pi N^j\nab^2\nab_j h)||_{\llb} $ can be bounded similarly.

 As for (\ref{tensor special}), since $r\leq 4$, $|\nab^sD_t^kv|$ can be bounded by
\begin{align}
|\nab^{s+1}D_t^{k-1}h|+C(K,M)(\sum_{j\leq r-1}|\nab^jv|+\sum_{k+s\leq r-1, s>0}|\p^sD_t^k h|).
\end{align}
and so (\ref{tensor special}) can be controlled via a priori assumptions times energies.

\prop \label{a priori estimate strong}
Let $r\geq r_0>\frac{n}{2}+\frac{3}{2}$, there is a continuous function $\mathcal{T}_r>0$ such that if
$$
0<T\leq \mathcal{T}_r(c_0,K,\mathcal{E}(0),E_r^*(0),\vol\Omega),
$$
where
\begin{align}
\mathcal{E}(t)=|(\nab_N h(t,\cdot))^{-1}|_{\linf}. \label{control ep}
\end{align}
Then any smooth solution of (\ref{EE}) for $0\leq t\leq T$ satisfies
\begin{align}
E_r^*(t)\leq 2E_r^*(0) \label{after gronwall},\\
\mathcal{E}(t)\leq 2\mathcal{E}(0)\label{check E 2},\\
g_{ab}(0,y)Z^aZ^b\lesssim g_{ab}(t,y)Z^aZ^b\lesssim g_{ab}(0,y)Z^aZ^b \label{eqv metric 2},
\end{align}
and there is some fixed $\eta>0$ such that
\begin{align}
|N(x(t,\bar{y}))-N(x(0,\bar{y}))| \lesssim \eta,\q\,\,\bar{y}\in \p\Omega \label{cont N},\\
|x(t,y)-x(0,y)| \lesssim \eta,\q\,\,y\in \Omega \label{cont x},\\
|\frac{\p x(t,\bar{y})}{\p y}-\frac{\p x(0,\bar{y})}{\p y}|\lesssim \eta,\q\,\,\bar{y}\in \p\Omega \label{cont dx}
\end{align}
hold.

To prove Proposition \ref{a priori estimate strong}, we will be using Sobolev lemmas. But then we must make sure that we can control the Sobolev constants. By Lemma \ref{interior sobolev} and \ref{boundary soboolev}, the Sobolev constants depend on $K = \frac{1}{l_0}$, in fact we are allowed to pick a $K$ depending only on initial conditions, which is proved in \cite{CL}. On the other hand, the change of the Sobolev constants in time are controlled by a bound for the time derivative of the metric in Lagrangian coordinate. We also need to control the constant $\frac{1}{\epsilon}$ appears to be in the Taylor sign condition (\ref{RT sign}).

\lem \label{check a priori} Assume the conditions in Proposition \ref{a priori estimate strong} hold. Then there are continuous functions $G_{r_0},H_{r_0},I_{r_0},J_{r_0}$ and $C_{r_0}$ such that
\begin{align}
||\nab v||_{L^{\infty}(\Omega)}\leq G_{r_0}(c_0,K,E_0,\cdots,E_{r_0}),\label{check 1}\\
||\nab h||_{L^{\infty}(\Omega)}\leq H_{r_0}(c_0,K,E_0,\cdots,E_{r_0},\vol\Omega),\label{check 2}\\
||\nab^2 h||_{L^{\infty}(\Omega)}+||\nab D_t h||_{L^{\infty}(\Omega)}\leq I_{r_0}(c_0,K, E_0,\cdots,E_{r_0},\vol\Omega),\label{check 2-deriv-h} \\
||\theta||_{\linf}\leq J_{r_0}(c_0, K,\mathcal{E},E_0,\cdots,E_{r_0},\vol\Omega) ,\label{check theta}\\
|\frac{d}{dt}\mathcal{E}|\leq C_{r_0}(c_0,K,\mathcal{E},E_0,\cdots,E_{r_0},\vol\Omega). \label{check E}
\end{align}

\begin{proof}
By Sobolev lemmas, we have
$$
||\nab v||_{L^{\infty}(\Omega)}\lesssim_K \sum_{1\leq j\leq 3} ||\nab^j v||_{\lli},
$$
$$
||\nab h||_{L^{\infty}(\Omega)}\lesssim_K \sum_{1\leq j\leq 3} ||\nab^j h||_{\lli},
$$
and
$$
||\nab^2 h||_{L^{\infty}(\Omega)}+||\nab D_t h||_{L^{\infty}(\Omega)}\lesssim_K \sum_{2\leq j\leq 4} ||\nab^j h||_{\lli}+\sum_{1\leq j\leq 3}||\nab^j D_th||_{\lli}.
$$
So, as a consequence of our interior and boundary estimates, (\ref{check 1})-(\ref{check 2-deriv-h}) follows. On the other hand, since $|\nab^2 h|\geq |\Pi \nab^2 h| = |\nab_N h||\theta|\geq \mathcal{E}^{-1}|\theta| $, so (\ref{check theta}) follows from (\ref{check 2-deriv-h}). Lastly, (\ref{check E}) is a consequence of (\ref{check 2-deriv-h}) and
$$
\frac{d}{dt}||(-\nab_N h(t,\cdot))^{-1}||_{\linf}\lesssim ||(-\nab_N h(t,\cdot))^{-1}||_{\linf}^2||\nab_N h_t(t,\cdot)||_{\linf} \label{change of e},
$$
\end{proof}

\subsubsection{Proof of Proposition \ref{a priori estimate strong}}

\indent Since when $r\geq r_0>\frac{n}{2}+\frac{3}{2}$, we have
$$
|\frac{d}{dt}E_r|\leq C_r(c_0,K,\mathcal{E},E_0,\cdots,E_{r_0},\vol \Omega)E_r^*,
$$
and the RHS is in fact a polynomial of $E_r^*$ with positive coefficients, we get (\ref{after gronwall}) from Lemma \ref{check a priori} and Gronwall's lemma if $\mathcal{T}_r(K,\mathcal{E}_0,E_r^*(0),\vol\Omega)>0$ is sufficiently small. (\ref{check E 2}) is a direct consequence of (\ref{check E}). In addition, we get from (\ref{after gronwall}) and Lemma \ref{check a priori} that
 \begin{align}
||\nab v||_{L^{\infty}(\Omega)}+||\nab h||_{L^{\infty}(\Omega)}\leq C(c_0,K,\mathcal{E}(0),E_0(0),\cdots,E_{r_0}(0),\vol \Omega), \label{nab v and h at 0}\\
||\nab^2 h||_{L^{\infty}(\Omega)}+||\nab D_th||_{L^{\infty}(\Omega)}+||\theta||_{\linf}\leq C(c_0,K,\mathcal{E}(0),E_0(0),\cdots,E_{r_0}(0),\vol \Omega). \label{nab^2 h and theta at 0}
\end{align}
It follows from these that, when $0<T\leq
\mathcal{T}_r(c_0,K,\mathcal{E}(0),E_r^*(0),\vol\Omega)$ with $\mathcal{T}_r$ chosen to be sufficiently small,
\begin{align}
||\nab v(t,\cdot)||_{L^{\infty}(\Omega)}+||\nab h(t,\cdot)||_{L^{\infty}(\Omega)}\lesssim ||\nab v(0,\cdot)||_{L^{\infty}(\Omega)}+||\nab h(0,\cdot)||_{L^{\infty}(\Omega)},
\end{align}
and
\begin{multline}
||\nab^2 h(t,\cdot)||_{L^{\infty}(\Omega)}+||\nab D_th(t,\cdot)||_{L^{\infty}(\Omega)}+||\theta(t,\cdot)||_{\linf}\\
\lesssim||\nab^2 h(0,\cdot)||_{L^{\infty}(\Omega)}+||\nab D_th(0,\cdot)||_{L^{\infty}(\Omega)}+||\theta(0,\cdot)||_{\linf},
\end{multline}
where $0<t\leq T$.

On the other hand, we have
\begin{align}
||v(t,\cdot)||_{L^{\infty}(\Omega)}\lesssim ||v(0,\cdot)||_{L^{\infty}(\Omega)}, \label{control v in the interior}\\
||\rho(t,\cdot)||_{L^{\infty}(\Omega)} \lesssim ||\rho(0,\cdot)||_{L^{\infty}(\Omega)}.\label{control h}
\end{align}
In fact, (\ref{control v in the interior}) follows since $D_tv=-\p h$ and (\ref{nab v and h at 0}), whereas (\ref{control h}) follows since $|D_t\rho|\leq |\rho\,\di v|$. Now,
  (\ref{eqv metric 2}) follows because $D_tg$ behaves like $\nab v$. Furthermore, (\ref{cont N}) follows from
$$
D_t N_a = -\frac{1}{2}N_a(D_t g^{cd})N_cN_d,
$$
and (\ref{eqv metric 2}). On the other hand, since by the definition of the Lagrangian coordinate, we have
$$
D_t x(t,y) = v(t,x(t,y)),
$$
and so (\ref{cont x}) follows since (\ref{control v in the interior}). Lastly, because
$$
D_t\frac{\p x}{\p y}=\frac{\p v(t,x)}{\p x}\frac{\p x}{\p y},
$$
(\ref{cont dx}) follows since (\ref{check 1}). \\
We close this section by briefly going over the idea which shows that one can choose $K$ depends only on the initial conditions.

\lem \label{l0 l1} Let $0\!\leq\! \eta\!\leq \!2$ be a fixed number, define $l_1\!=\!l_1(\eta)$ to be the largest number such that
$$
|N(\bar{x}_1)-N(\bar{x}_2)|\leq \eta, \q\text{whenever}\,\, |\bar{x}_1-\bar{x}_2|\leq l_1, \bar{x}_1,\bar{x}_2\in\p\DD_t.
$$
Suppose $|\theta|\leq K$, we recall that $l_0$ is the injective radius defined in Section 1.3, then
\begin{align*}
l_0\geq \min(l_1/2,1/K),\\
l_1\geq \min(2l_0, \eta/K).
\end{align*}
\begin{proof}
See Lemma 3.6 of \cite{CL}
\end{proof}
In fact, Lemma \ref{l0 l1} shows that $l_0$ and $l_1$ are comparable as long as the free surface is regular.
\lem \label{K initial data}
Fix $\eta>0$ sufficiently small, let $\mathcal{T}$ be in Proposition \ref{a priori estimate strong}. Pick $l_1>0$ such that, whenever $|x(0,y_1)-x(0,y_2)|\leq 2l_1$,
\begin{align}
|N(x(0,y_1))-N(x(0,y_2))|\leq \frac{\eta}{2} \label{K initial}.
\end{align}
Then if $t\leq \mathcal{T}$ we have
\begin{align}
|N(x(t,y_1))-N(x(t,y_2))|\leq \eta, \label{K at t}
\end{align}
whenever $|x(t,y_1)-x(t,y_2)|\leq l_1$.
\begin{proof}
We have
\begin{multline}
|N(x(t,y_1))-N(x(t,y_2))|\\
\leq |N(x(t,y_1))-N(x(0,y_1))|+|N(x(0,y_1))-N(x(0,y_2))|+|N(x(0,y_2))-N(x(t,y_2))|,
\end{multline}
and so (\ref{K at t}) follows from (\ref{cont N}) and (\ref{cont x}).
\end{proof}
Lemma \ref{K initial data} allows us to pick
$
l_1(t)\leq \frac{l_1(0)}{2},
$
in other words, we have if $\frac{1}{l_1(0)}\leq \frac{K}{2}$, then
$
\frac{1}{l_1(t)}\leq K.
$
Therefore, Lemma \ref{l0 l1} yields
$
\frac{1}{l_0(t)}\leq K.
$

\section{Incompressible limit}\label{section 6}
In this section we prove that the energy estimates for compressible Euler equations are in fact uniform in the sound speed. For physical reasons, the sound speed is defined by
$$
 c(t,x)^2:= p'(\rho)
$$
Let $\{p_{\kk}(\rho)\}$ be a family parametrized by $\kk\in\RR^+$, such that for each $\kk$ we have
$$
p_{\kk}'(\rho)|_{\rho=1}=\kk.
$$
We shall call $\kk$ as the sound speed by a slight abuse of terminology. We are concerning with fluid motion when $\kk$ is large and in its limit as $\kk\to\infty$.
We recall that the the enthalpy $h$ has derivative
$$
h'(\rho) = \frac{p_{\kk}'(\rho)}{\rho},
$$
and since $p_{\kk}(\rho)$ is strictly increasing for every $\kk$ and $h'(\rho)>0$, we can write $\rho$ as a function of $h$ depends on $\kk$. We want to impose the following conditions on $\rho_{\kk}(h)$:
\begin{equation}
\rho_{\kk}(h)\to 1, \quad\text{and}\quad e_{\kk}(h)\to 0,\qquad \text{as}\quad
\kk\to\infty
\end{equation}
and for some fixed constant $c_0$
\begin{equation}\label{econdkappa}
|e^{(k)}_{\kk}(h)|\leq c_0,\qquad\text{and}\qquad
|e^{(k)}_{\kk}(h)|\leq c_0\sqrt{e^{\,\prime}_{\kk}(h)}, \qquad\text{for $k\leq 6$}.
\end{equation}

The purpose of this section is to prove:
\prop \label{uniform kk}
Let $E_{r,\kk}$ be defined as (\ref{Er}) and for $r\geq r_0>\frac{n}{2}+\frac{3}{2}$, there is a continuous function $\mathcal{T}_r>0$ such that if
$$
0<T\leq \mathcal{T}_r(c_0,K,\frac{1}{\epsilon},E_{r,\kk}^*(0),\vol\Omega),
$$
then any smooth solution of (\ref{EE}) for $0\leq t\leq T$ satisfies
\begin{align}
E_{r,\kk}^*(t)\leq 2E_{r,\kk}^*(0),
\end{align}
provided the physical sign condition
$$
-\nab_{N}h_{\kk}\geq \epsilon>0,\q\text{on}\,\p\Omega
$$
holds.
\rmk $\mathcal{T}_\kk$ depends on $\kk$ only through $E_{r,\kk}^*(0)$, and in fact we show in Section \ref{section 7} that there is a sequence of data $(v_{0,\kk}, h_{0,\kk})$ for which the initial energy $E_{r,\kk}^*(0)$ is uniformly bounded independent of $\kk$, as $\kk\to \infty$.

\rmk Our energy estimate can be slightly modified so that it can be carried over to the incompressible Euler equations (i.e., the case when $\kk =\infty$). We consider
\begin{align}
\widetilde{E}_r = \sum_{s+k=r}E_{s,k}+K_r+\widetilde{W}_{r+1},\label{tilde E_r}
\end{align}
where $E_{s,k}$ and $K_r$ are defined as in \eqref{Esk}-\eqref{K_r}, but
\begin{align}
\widetilde{W}_{r+1} = \frac{1}{2}||e'(h)D_t^{r+1}h||_{\lli}+\frac{1}{2} ||\sqrt{e'(h)}\nab D_t^{r}h||_{\lli}.
\end{align}
In addition, we assume
$$
|e_{\kk}^{(k)}(h)|\leq c_0|e_{\kk}'(h)|, \q\text{for}\,\, k\leq 6.
$$
The energy estimates for \eqref{tilde E_r} follows from Theorem \ref{the key lemma}, with the norms $||h||_r$ and $||D_th||_r$ being replaced by a weaker version, i.e.,
$$
||h||_{r, wk} = \sum_{s+k=r, k\leq r-2}||\nab^s D_t^k h||_{\lli}+||e'(h)D_t^{r}h||_{\lli}+ ||\sqrt{e'(h)}\nab D_t^{r-1}h||_{\lli},
$$
and thus
$$
||D_th||_{r, wk} = \sum_{s+k=r, k\leq r-2}||\nab^s D_t^{k+1} h||_{\lli}+ ||e'(h)D_t^{r+1}h||_{\lli}+ ||\sqrt{e'(h)}\nab D_t^{r}h||_{\lli}.
$$
In addition to this, in view of Section \ref{section 5.2}, we need to estimate $||\nab D_t^{j}h||_{\lli}, 0\leq j\leq r-1$ independently, without using $W_{r+1}$. This follows from
\begin{align*}
||\nab D_t^j h||_{\lli} \lesssim_{\vol\Omega} ||e'(h)D_t^{j+2}h||_{\lli}+||f_{j+1}||_{\lli}+||g_{j+1}||_{\lli}.
\end{align*}
The first term is part of $\widetilde{W}_{r+1}^*$, whereas the bounds for $||f_{j+1}||$ and $||g_{j+1}||$ do not rely on $||\nab D_t^j h||_{\lli}$ (see Section \ref{section 5.5}),  and so they can be bounded by $||D_th||_{r, wk}$ (plus lower order terms) and thus $\widetilde{E}_r^*$ via the estimates in Section \ref{section 5.3}.

\rmk Related the issue above, in the incompressible case,
i.e. when $e(h)=0$, one can actually bound $W_{r+1}(0)$ by
$E_{s,k}(0)$ for $s+k\leq r$ and $K_r(0)$ provided that $W_{r+1}(0)<\infty$. In fact, since then
$\triangle h=-\p v\cdot \p v$, it follows that to highest order
$\triangle D_t^r h=-\p v\cdot \p D_t^r v\sim -\p v\cdot \p^2 D_t^{r-1} h$.
Multiplying both sides by $D_t^r h $ and integrating by parts the second derivatives on both sides gives that
to highest order
$||\nab D_t^r h ||_{\lli}^2\sim ||\nab D_t^r h ||_{\lli} ||\nab D_t^{r-1} h ||_{\lli}$. If $||\nab D_t^r h ||_{\lli}<\infty$, as will be the case for smooth data or limits of smooth data, we conclude that
$||\nab D_t^r h ||_{\lli}$ can in fact be bounded by $ ||\nab D_t^{r-1} h ||_{\lli}$
and lower order terms. This shows that it makes sense to include $W_{r+1}$ in the energy also in the incompressible case if one needs to keep control of time
derivatives to highest order.
\\

Based on the analysis we have in Section \ref{section 5.5}, Proposition \ref{uniform kk} is a direct consequence of:
\thm\label{uniform energybounds}
Let $E_{r,\kk}$ be defined as (\ref{Er}), then there are continuous functions $C_r$ such that, for $t\in[0,T]$ and $r\leq 4$ that
\begin{align}
|\frac{dE_{r,\kk}(t)}{dt}|\leq C_r(K,\frac{1}{\epsilon},M,c_0,\vol\DD_t,E_{r-1,\kk}^*)E_{r,\kk}^*(t) \label{uniform kk energy est},
\end{align}
for all $\kk$, provided that the assumptions on $e_{\kk}(h)$ hold and
\begin{align}
|\theta_{\kk}|+\frac{1}{l_0}\leq K,\q\text{on}\,\p\Omega,\label{boot kk 1}\\
-\nab_{N}h_{\kk}\geq \epsilon>0,\q\text{on}\,\p\Omega,\label{RE kk}\\
1\leq \rho_{\kk}\leq M,\q\text{in}\,\Omega,\\
|D_t^2h_{\kk}|\leq M,\q\text{in}\,\Omega,\label{boot kk extra}\\
|\nab v_{\kk}|+|\nab h_{\kk}|+|\nab^2 h_{\kk}|+|\nab D_th_{\kk}|\leq M,\q\text{in}\,\Omega. \label{boot kk v,h}
\end{align}
\rmk We actually do not need to assume the bound for $|D_th_{\kk}|$. Since $\Omega$ is bounded and $D_th_{\kk}=0$ on $\p\Omega$, so
$$
|D_th_{\kk}|_{L^{\infty}} \lesssim \int_{\Omega} |\nab D_th_{\kk}| \lesssim_{\vol\Omega} M.
$$
Together with (\ref{boot kk extra}), we have
\begin{align}
|D_th_{\kk}|+|D_t^2h_{\kk}|\leq M,\q\text{in}\,\Omega,
\end{align}
independent of $\kk$. This is compatible with the case with fixed sound speed. \\

\indent To prove Proposition \ref{uniform kk}, the analysis in the Section 5 implies that it suffices to prove that $||v_{\kk}||_{r,0}$, $||h_{\kk}||_r$ and $\lee h_{\kk}\ree_r$ are bounded uniformly in $\kk$. It is easy to see that under a priori assumptions (\ref{boot kk 1})-(\ref{boot kk v,h}), the estimates for $||f_r||_{\lli}$ and $||g_r||_{\lli}$ (Section 4) stay unchanged.

The interior estimates in Section 5.2 are uniform in the sound speed since $\sum_{i\leq r}||h_{\kk}||_i$ involves terms of the form $\sum_{i\leq r}||\sqrt{e_{\kk}'(h)}D_t^ih_{\kk}||_{\lli}$ for each $r$, which means that we do not need the lower bound of $|e_{\kk}'(h)|$ in our estimates. Further, the boundary estimates for $\sum_{i\leq r}\lee h_{\kk}\ree_i$ follows as well, which are uniform in $\kk$ since the interior estimates are. But we need (\ref{g_r mod 2}) to estimate $\lee h_{\kk}\ree_4$. Finally, as for the extra a priori assumption (\ref{boot kk extra}), we can get it back by the interior estimates as in Lemma \ref{check a priori} since
\begin{align}
||D_t^2h_{\kk}||_{L^{\infty}(\Omega)}\lesssim_K \sum_{j=1,2}||\nab^j D_t^2h_{\kk}||_{\lli}+||D_t^2h_{\kk}||_{\lli},
\end{align}
where $||D_t^2h_{\kk}||_{\lli}\lesssim_{\vol\Omega} ||\nab D_t^2h_{\kk}||_{\lli}$ via (\ref{FB}).

The above analysis shows that
$$
E_{r,\kk}^*(t)\leq 2E_{r,\kk}^*(0),
$$
regardless of the sound speed $\kk$. Furthermore, since we are able to show that there exists a sequence of data $(v_{0,\kk},h_{0,\kk})$ such that $E_{r,\kk}^*(0)$ are uniformly bounded in Section \ref{section 7}. A direct consequence of this is that $v_{\kk}$ and $h_{\kk}$ converge in $C^2([0,T],\Omega)$ as $\kk\to \infty$. To be more precise, we define
\mydef
The space
$C^l([0,T],\Omega)$
consists all functions $u(t,x)$ with
$\nab^s D_t^k u(t,\cdot),\,\,s+k\leq l$
continuous in $\Omega$.

\thm \label{main 2}
 Let $u_0$ be a divergence free vector field such that its corresponding pressure $p_0$, defined by $\lap p_0 = -(\p_i u_0^k)(\p_k u_0^i)$ and $p_0\big|_{\p\DD_0}=0$, satisfies the physical condition $-\nab_N p_0\big|_{\p\DD_0} \geq \epsilon >0$. Let $(u, p)$ be the solution of the incompressible free boundary Euler equations with data $u_0$, i.e.
$$
\rho_0 D_t u = -\p p,\q \di u=0,\qquad p|_{\p\DD_0}=0,\quad u|_{t=0}=u_0
$$
with the constant density $\rho_0 =1$. Furthermore, let $(v_{\kk}, h_{\kk})$ be the solution for the compressible Euler equations \eqref{EEkk}, with the density function $\rho_{\kk} :h\rightarrow \rho_{\kk}(h)$, and the initial data $v_{0\kk}$ and $ h_{0\kk}$, satisfying the compatibility condition up to order $r+1$, as well as the physical sign condition \eqref{RE kk}. Suppose that $\rho_{\kk} \to \rho_0=1$, $v_{0\kk}\to u_0$ and $h_{0\kk}\to p_0$ as $\kk\to\infty$, such that $E_{r,\kk}^*(0)$ is bounded uniformly independent of $\kk$, then
$$
(v_{\kk},h_{\kk})\to (u,p)\q \text{in}\,\,C^2([0,T],\Omega).
$$
\begin{proof}
We first show that the $C^2$ norms of $v_{\kk}$ and $h_{\kk}$ are bounded by $E_{r,\kk}^*(t)$. By the definition of $C^2([0,T],\Omega)$, we have
\begin{align*}
||v_{\kk}||_{C^2([0,T],\Omega)}^2+||h_{\kk}||_{C^2([0,T],\Omega)}^2 \lesssim_K \sum_{j\leq 4}(||v_{\kk}||_{j,0}^2+||h_{\kk}||_{j,0}^2)\lesssim E_{4,\kk}^*(t)\leq 2E_{4,\kk}^*(0),
\end{align*}
by Sobolev lemma. Hence, the energy estimates (\ref{uniform kk energy est}) as well as the arguments in Section \ref{section 7.3} yield that the quantities
$||v_{\kk}||_{C^2([0,T],\Omega)}$ and $||h_{\kk}||_{C^2([0,T],\Omega)}$
are uniformly bounded. Furthermore, since the uniform bound for $\sum_{j\leq 4}(||v_{\kk}||_{j,0}+||h_{\kk}||_{j,0})$ also implies that for $s+k\leq 2$
$$
\nab^s D_t^k v_{\kk},\,\, \nab^s D_t^k h_{\kk}\in C^{0,\frac{1}{2}}(\Omega),
$$
in the sense of H\"{o}lder continuous functions (see \cite{Ev}, Chapter 5). This shows that the families $v_{\kk}$ and $h_{\kk}$ are in fact equi-continuous in $C^2([0,T]\times \Omega)$.
Therefore, $(v_{\kk},h_{\kk})$ convergent (after possibly passing to subsequence) by Arzela-Ascoil theorem. Finally, $(v_{\kk},h_{\kk})\to (u, p)$ follows from $e'(h)D_t^l h_{\kk}\to 0$ for $l=1,2$, which is a direct consequence of $||h_{\kk}||_{C^2([0,T],\Omega)}^2$ is bounded uniformly.

\end{proof}

\section{Existence of initial data satisfying the compatibility condition}\label{section 7}

In this section we show that given any incompressible data there is a sequence of compressible initial data, depending on the sound speed $\kk$, that satisfy the compatibility conditions and converges to the given incompressible data in our energy norm, as the sound speed $\kk\to\infty$. Hence by the previous theorem the incompressible limit will exist for this sequence.

Given $u_0$ a divergence free vector field such that its corresponding pressure $p_0$, defined by $\lap p_0 = -(\p_i u_0^k)(\p_k u_0^i)$ and $p_0\big|_{\partial\Omega}=0$, satisfies the physical condition $-\nab_N p_0\big|_{\partial \Omega} \geq \epsilon >0$, we are going to construct a sequence of incompressible data $(v_{0},h_{0})=(v_{0\kk},h_{0\kk})$ satisfying the compatibility conditions such that the corresponding solutions converge to the solution of the incompressible equations with data $(u_0,p_0)$ in the energy norm
initially, as the sound speed $\kk\to \infty$.

For simplicity we assume that $e(h) =\kk^{-1}h$. We consider the compressible Euler's equations
\begin{align}
D_t v&=-\partial h,\label{eulersec7}\\
\kk^{-1} D_t h& =-\di v,\label{continuitysec7}
\end{align}
in $\Omega$ (in the Lagrangian coordinates) with boundary condition
\begin{equation}
 h\big|_{\partial\Omega}=0, \label{boundaryconditionsec7}
\end{equation}
and initial data
\begin{equation}
v\big|_{t=0}=v_0,\qquad h\big|_{t=0}=h_0, \label{initialdatasec7}
\end{equation}
depending on $\kk$.
In order for initial data to be compatible with the boundary condition we must have
\begin{equation}
h_0\big|_{\partial\Omega}=0,\qquad \di v_0\big|_{\partial\Omega}=0,
\end{equation}
since we must also have that $D_t h\big|_{\partial\Omega}=0$ at time $0$.
Moreover since $h$ satisfies the wave equation
\begin{align}
\kk^{-1} D_t^2h = \lap h + (\p_i v^k)(\p_k v^i),
\end{align}
and $D_t^2 h\big|_{\partial \Omega}=0$ when $t=0$, we must also have
\begin{align}
\lap_0 h_0 + (\p_i v_0^k)(\p_k v_0^i)= 0,\q \text{on}\q \p\Omega. \label{lowest comp}
\end{align}
Here $\lap_0$ is the Laplacian with respect the smooth metric \eqref{g} at time $0$ on the domain with smooth boundary $\partial\Omega$, and $\partial_i =\partial y^a/\partial x^k \partial/\partial y^a$ is a smooth differential operator at time $0$.
Similarly, by differentiating the wave equation we get
\begin{align}
\kk^{-1} D_t^3h = \lap D_t h + f_{2},
\end{align}
for some $f_2$ as in section \ref{section 4}.
Since we also want $D_t^3 h\big|_{\partial \Omega}=0$ when $t=0$ we also need
\begin{align}
\lap_0 h_1 + F_1 = 0,\q \text{on}\q \p\Omega\label{next lowest comp},\qquad\text{where}\quad
h_1=D_t h\big|_{t=0}
\end{align}
 and $F_1 = f_2|_{t=0}$ is a function of $v_0,h_0$ and its space derivatives. Similarly we get
 \begin{align}
\kk^{-1} D_t^{k+2} h = \lap D_t^k h + f_{k+1},
\end{align}
 and hence we must have 
 \begin{align}
\lap_0 h_k + F_k = 0,\q \text{on}\q \p\Omega,\qquad\text{where}\quad
h_k=D_t^k h\big|_{t=0}
\end{align}
 and $F_k = f_{k+1}|_{t=0}$ if a function of $v_0, h_0,..., h_{k-1}$ and its space derivatives.

Given a divergence free vector field $u_0$, let
\begin{equation}
v_{0} = u_0+\p \phi.\label{u0}
\end{equation}
Then the continuity equation requires that
\begin{equation}
\lap_0\phi  = -\kk^{-1}h_1, \label{Laplacephi}
\end{equation}
and we will choose boundary conditions, e.g.
\begin{equation}
\nab_N\phi |_{\p\Omega}=0.\label{boundaryconditionphi}
\end{equation}
Moreover the time derivatives of the wave equation require that
\begin{equation}
\lap_0 h_k + F_k = \kk^{-1} h_{k+2},\quad \text{in}\quad \Omega \quad \text{and}\quad h_k|_{\p\Omega} =0,\qquad k=0,\dots,N\label{timedifferentiatewaveequation}
\end{equation}
where $F_k$ are function of $v_0, h_0,..., h_{k-1}$ and its space derivatives.
If we prescribe $h_{N+1}$ and $h_{N+2}$ to be any functions that vanish at the boundary,
e.g.
\begin{equation}
h_{N+1}=h_{N+2}=0, \qquad \text{in}\quad \Omega . \label{higestderivativeconditions}
\end{equation}
Then \eqref{u0}-\eqref{higestderivativeconditions}
gives a system for $(v_0, h_0, h_1,\cdots, h_N,h_{N+1},h_{N+2})$, such that  when
$\kk\to \infty $ the compressible data $(v_0,h_0)\to (u_0,p_0)$, the incompressible data,
and for each $\kk$, $(v_0,h_0)$ satisfy the $N$ compatibility conditions. It remains to show that the system \eqref{u0}-\eqref{higestderivativeconditions} has a solution if $\kk$ is sufficiently large with uniformly bounded energy norms as $\kk\to\infty$.

\subsection{The a priori energy bounds for the full system}
Our energy estimate requires that the compatibility conditions to be satisfied up to $5$th order, i.e., we
need to find $v_{0}$, and $D_t^kh_{\kk}|_{t=0} = h_k$ such that
\begin{align}
h_k|_{\p\Omega} = 0,\q 0\leq k\leq 5.
\end{align}
This can be achieved by solving \footnote{The Neumann boundary condition $\nab_N\phi_{\kk} |_{\p\Omega}=0$ can be replaced by the Dirichlet boundary condition $\phi_{\kk}|_{\p\Omega}=0$. Nevertheless, the Neumann condition makes more sense here since it does not change the boundary velocity. In addition, one may think $\phi$ as $h_{-1}$ and so that it would be more natural if we impose the Dirichlet boundary condition in view of this. }
\begin{align}
\begin{cases}
v_{0} = u_0+\p \phi,\q \text{in}\,\,\Omega,\\
\lap\phi = -\kk^{-1}h_1,\q \text{in}\,\,\Omega,\q\text{and}\q\nab_N \phi|_{\p\Omega}=0,\\
\lap h_k = F_k + \kk^{-1} h_{k+2},\q \text{in}\,\,\Omega,\q\text{and}\q h_k|_{\p\Omega}=0, \,\,0\leq k\leq 3,\\
h_4 = h_5=0,\q \text{in}\,\,\Omega.\label{3-system}
\end{cases}
\end{align}
 Here,
\begin{align}
F_k=c_{\alpha_1\cdots\alpha_m\beta_1\cdots\beta_n}^{\gamma_1\cdots\gamma_n,k}(\p^{\alpha_1} v_{0})\cdots(\p^{\alpha_m} v_{0})(\p^{\beta_1}h_{\gamma_1})\cdots(\p^{\beta_n}h_{\gamma_n}),\q k=2,3 \label{F_3}
\end{align}
where
\begin{align*}
\alpha_1+\cdots+\alpha_m+(\beta_1+\gamma_1)+\cdots+(\beta_n+\gamma_n) = k+2,\\
1\leq \alpha_i \leq k ,\q 1\leq \beta_j+\gamma_j \leq k+1,\q \beta_j\geq 1,\q 0\leq \gamma_j\leq k-1 ,\\
1\leq m+n\leq k+2.
\end{align*}
Here, either $m$ or $n$ can be $0$, in which case it means that the corresponding factor is not involved in the product. In addition, 
\begin{align}
F_1 = c_{\alpha,\beta,1}(\p^{\alpha}v_{0})(\p^{\beta}h_0)+d_1(\p v_{0})^3, \q 1\leq \alpha,\beta\leq 2, \q \alpha+\beta=3,\\
F_0= (\p v_{0})^2.
\end{align}

We show the existence of solution for (\ref{3-system}) by successive approximation starting from the solution $(h_0^0,h_1^0,h_2^0,h_3^0)$ that solves
\begin{align}
\lap h_k^0 = F_k(\p^{\alpha} v_0, \p^{\beta_0}h_0^0,\cdots,\p^{\beta_{k-1}}h_{k-1}^0),\q k=0,1,2,3 \label{3-system nu=0}
\end{align}
and we define $(v_{0}^\nu,h_0^{\nu},\cdots,h_3^{\nu})$ inductively by solving
\begin{align}
\begin{cases}
v_{0}^\nu = u_0 +\p\phi^\nu,\\
\lap \phi^\nu = -\kk^{-1}h_1^{\nu-1},\\
\lap h_0^{\nu} = \kk^{-1} h_2^{\nu-1}+(\p v_{0}^\nu)(\p v_{0}^\nu),\\
\lap h_1^{\nu} = \kk^{-1} h_3^{\nu-1}+F_1(\p ^{\alpha}v_{0}^\nu,\p^{\beta_0}h_0^{\nu}),\\
\lap h_2^{\nu} = F_2(\p^{\alpha} v_{0}^\nu, \p^{\beta_0}h_0^{\nu},\p^{\beta_1}h_1^{\nu}),\\
\lap h_3^{\nu} = F_3(\p^{\alpha} v_{0}^\nu, \p^{\beta_0}h_0^{\nu},\p^{\beta_1}h_1^{\nu},\p^{\beta_2}h_2^{\nu}),\\
 h_k^\nu|_{\p\Omega} = \nab_N \phi^\nu|_{\p\Omega}=0\q \text{for all}\,\,\nu,\label{3-system nu}
\end{cases}
\end{align}
for $\nu\geq 1$.

We define that for $0\leq k\leq 3$,
\begin{align*}
m_k^{\nu} := ||h_k^{\nu}||_{H^{s-k}(\Omega)},\q s\geq 5\\
m_*^{\nu} := \sum_k m_k^{\nu}+||v_0^{\nu}||_{H^s}.
\end{align*}
We shall apply the standard elliptic estimate as well as Sobolev lemmas to get bounds for $m_k^{\nu}$. It is worth mentioning that since we are working in dimensions $2$ and $3$, $H^l$ is an algebra for $l\geq 2$, i.e.,
\begin{align}
||fg||_{H^l} \leq C||f||_{H^l}||g||_{H^l} \label{sobolev estimates},
\end{align}
which is a direct consquence of Sobolev lemma (\ref{interior sobolev 2s>n}). On the other hand, when $l<2$, Sobolev lemma (\ref{interior sobolev 2s<n}) implies
\begin{align}
||u_1\cdots u_N||_{H^l} \leq
\begin{cases}
 C||u_1||_{H^{l+1}}\cdots||u_N||_{H^{l+1}},\q N=2,3 \\
  C||u_1||_{H^{l+2}}\cdots||u_N||_{H^{l+2}},\q N\geq 4 \label{sobolev estimates'}.
\end{cases}
\end{align}

Now, since $s\geq 5$ and $\Omega\in \RR^n$, $n=2,3$ imply that each $H^{l}, l\geq s-3$ is an algebra, we have by \eqref{3-system nu} that
\begin{align}
m_0^{\nu} \leq C\kk^{-1}||h_2^{\nu-1}||_{H^{s-2}}+C||(\p v_{0}^\nu)^2||_{H^{s-2}} \leq C(\kk^{-1} m_2^{\nu-1}+||v_{0}^\nu||_{H^{s}}^2) \label{m_0^nu},
\end{align}
and
\begin{dmath}
m_1^{\nu} \leq C\kk^{-1} ||h_3^{\nu-1}||_{H^{s-3}}+C||F_1||_{H^{s-3}} \leq C(\kk^{-1} m_3^{\nu-1}+||v_{0}^\nu ||_{H^{s}}||h_0||_{H^s}+||v_{0}^\nu ||_{H^{s}}^3) = C(\kk^{-1} m_3^{\nu-1}+||v_{0}^\nu ||_{H^{s}}m_0^{\nu}+||v_{0}^\nu ||_{H^{s}}^3)\label{m_1^nu}.
\end{dmath}
However, the estimates for $m_2^{\nu}$ and $m_3^{\nu}$ are a bit more complicated. The standard elliptic estimate yields that
\begin{align*}
m_2^{\nu}\leq C||F_2||_{H^{s-4}},\\
m_3^{\nu} \leq C||F_3||_{H^{s-5}},
\end{align*}
and the following analysis is devoted to bound $||F_2||_{H^{s-4}}$ and $||F_3||_{H^{s-5}}$.
\subsubsection{Bounds for $||F_2||_{H^{s-4}}$}
Since $F_2$ is a sum of products of the form \eqref{F_3} with $k=2$, and each product involves at least $2$ but no more than $4$ terms, we have
\begin{itemize}
\item If the product involves less than 4 terms, i.e., it is of the form
\begin{align*}
(\p^{\alpha_1} v_{0})\cdots(\p^{\alpha_m} v_{0})(\p^{\beta_1}h_{\gamma_1})\cdots(\p^{\beta_n}h_{\gamma_n}),\q m+n\leq 3,
\end{align*}
and for each $i,j$,  we have $1\leq \alpha_i\leq 2$, $1\leq (\beta_j+\gamma_j)\leq 3$
as well as
\begin{align*}
\alpha_1+\cdots+\alpha_n+(\beta_1+\gamma_1)+\cdots+(\beta_n+\gamma_n) = 4.
\end{align*}
This guarantees that only
\begin{align*}
\p^{\alpha} v_{0},\q 1\leq \alpha\leq 2,\\
\p^{\beta_k}h_k,\q 1\leq \beta_k\leq 3-k,\q k=0,1
\end{align*}
are allowed in the product.
Therefore, by (\ref{sobolev estimates'}), we get
\begin{dmath}
||(\p^{\alpha_1} v_{0})\cdots(\p^{\alpha_m} v_{0})(\p^{\beta_1}h_{\gamma_1})\cdots(\p^{\beta_n}h_{\gamma_n})||_{H^{s-4}}\\
 \leq C||\p^{\alpha_1}v_{0}||_{H^{s-3}}\cdots ||\p^{\alpha_m}v_{0}||_{H^{s-3}}||\p^{\beta_1}h_{\gamma_1}||_{H^{s-3}}\cdots||\p^{\beta_n}h_{\gamma_n}||_{H^{s-3}}\\
 \leq Cp_2(||v_{0}||_{H^{s}}, ||h_0||_{H^s}, ||h_1||_{H^{s-1}}),
\end{dmath}
for some polynomial $p_2$.
\item If the product involves exactly $4$ terms, then it must be $(\p v_0)^4$.  Hence, by the Sobolev estimate (\ref{sobolev estimates}),
\begin{dmath}
||(\p v_0)^4||_{H^{s-4}}\leq ||(\p v_0)^4||_{H^{s-3}} \leq ||v_0||_{H^s}^4.
\end{dmath}
\end{itemize}
Therefore, we conclude that
\begin{align}
m_2^{\nu} \leq CP_2( ||v_{0}^\nu||_{H^{s}}, m_0^{\nu},m_1^{\nu}),
\end{align}
for some polynomial $P_2$.
\subsubsection{Bounds for $||F_3||_{H^{s-5}}$}
We shall proceed as in the previous case. Since $F_3$ is a sum of products of the form (\ref{F_3}) with $k=3$, and each product involves at least $2$ but no more than $5$ terms, we have
\begin{itemize}
\item If the product involves no more than $3$ terms, i.e., it is of the form
\begin{align*}
(\p^{\alpha_1} v_{0})\cdots(\p^{\alpha_m} v_{0})(\p^{\beta_1}h_{\gamma_1})\cdots(\p^{\beta_n}h_{\gamma_n}),\q 1\leq m+n\leq 3,
\end{align*}
and for each $i,j$,  we have $1\leq \alpha_i\leq 3$, $1\leq (\beta_j+\gamma_j)\leq 4$
as well as
\begin{align*}
\alpha_1+\cdots+\alpha_n+(\beta_1+\gamma_1)+\cdots+(\beta_n+\gamma_n) = 5.
\end{align*}
This implies that only
\begin{align*}
\p^{\alpha} v_{0},\q 1\leq \alpha\leq 3,\\
 \p^{\beta_k} h_k,\q 1\leq \beta_k \leq 4-k,\q k=0,1,2
\end{align*}
are allowed to be included in the product. Therefore, by (\ref{sobolev estimates'}), we have
\begin{dmath}
||(\p^{\alpha_1} v_{0})\cdots(\p^{\alpha_m} v_{0})(\p^{\beta_1}h_{\gamma_1})\cdots(\p^{\beta_n}h_{\gamma_n})||_{H^{s-5}}\\
 \leq C||\p^{\alpha_1}v_{0}||_{H^{s-4}}\cdots ||\p^{\alpha_m}v_{0}||_{H^{s-4}}||\p^{\beta_1}h_{\gamma_1}||_{H^{s-4}}\cdots||\p^{\beta_n}h_{\gamma_n}||_{H^{s-4}}\\
 \leq Cp_3(||v_{0}||_{H^{s}},||h_0||_{H^s},||h_1||_{H^{s-1}},||h_2||_{H^{s-2}}),
\end{dmath}
for some polynomial $p_3$.
\item If the product involves exactly $4$ terms, then for each $i,j$,
\begin{align*}
1\leq \alpha_i \leq 2,\\
1\leq \beta_j+\gamma_j\leq 2,
\end{align*}
and thus only
\begin{align*}
\p^{\alpha} v_{0},\q 1\leq \alpha\leq 2,\\
\p^{\beta_k} h_k,\q 1\leq \beta_k\leq 2-k,\q k=0,1\\
\end{align*}
are allowed in each product. Hence one can then use Sobolev estimate (\ref{sobolev estimates}) to bound the product in $H^{s-3}$ norm, i.e.,
\begin{align}
||(\p^{\alpha_1} v_{0})\cdots(\p^{\alpha_m} v_{0})(\p^{\beta_1}h_{\gamma_1})\cdots(\p^{\beta_n}h_{\gamma_n})||_{H^{s-3}}\leq C p_3(||v_{0}||_{H^{s-1}}, ||h_0||_{H^s}, ||h_1||_{H^{s-1}}).
\end{align}
\item If the product involves exactly $5$ terms, then it must be $(\p v_0)^5$.  Hence, by the Sobolev estimate (\ref{sobolev estimates}),
\begin{dmath}
||(\p v_0)^5||_{H^{s-5}}\leq ||(\p v_0)^5||_{H^{s-3}} \leq ||v_0||_{H^s}^5.
\end{dmath}
\end{itemize}
Therefore, we conclude that
\begin{align}
m_3^{\nu} \leq CP_3(||v_{0}^\nu||_{H^{s}},m_0^{\nu},m_1^{\nu},m_2^{\nu}),
\end{align}
for some polynomial $P_3$.

Summing up the estimates for $m_k^{\nu}$ for $0\leq k\leq 3$, we have
\begin{align}
\sum_k m_k^{\nu} \leq P(C, \kk^{-1} \sum_k m_k^{\nu-1}, ||v_{0}^\nu||_{H^{s}})\label{sum mnu},
\end{align}
for some polynomial $P$.  In particular, this implies that $m_*^0=\sum_k m_k^0$ is uniformly bounded by a function depends on $||u_0||_{H^{s}}$ by taking $\nu=0$. Moreover, since $\lap\phi_{\kk}^{\nu} = -\kk h_1^{\nu-1}$ and $\nab_N\phi_{\kk}^{\nu}|_{\p\Omega}=0$, we have
\begin{align}
|| \p\phi_{\kk}^{\nu}||_{H^s} \leq C\kk^{-1}||h_1^{\nu-1}||_{H^{s-1}} = C\kk^{-1}m_1^{\nu-1}.
\end{align}
Hence,
\begin{align}
||v_{0}^{\nu}||_{H^{s}} \leq C(||u_0||_{H^{s}}+||\p \phi_{\kk}^\nu||_{H^{s}})\leq C(||u_0||_{H^s}+\kk^{-1}m_1^{\nu-1}).
\end{align}

This, together with \eqref{sum mnu} yield
\begin{align}
\sum_k m_k^\nu+||v_0^\nu||_{H^s}=m_*^{\nu} \leq P_1(C, \kk^{-1} m_*^{\nu-1}, ||u_{0}||_{H^{s}}),
\end{align}
and so it is uniformly bounded for all $\nu$ by induction if $\kk^{-1}$ is chosen to be sufficiently small.

\subsection{The iteration scheme  }
Let's define
\begin{align*}
V^\nu := v_{0}^\nu-v_{0}^{\nu-1},\\
\Phi^\nu:=\phi^\nu-\phi^{\nu-1},\\
A_k^{\nu} := h_k^{\nu} - h_k^{\nu-1},\\
M_k^{\nu} := ||A_k^{\nu}||_{H^{s-k}(\Omega)},\q s\geq 5,\\
M_*^{\nu} := \sum_k M_k^{\nu}+||V^\nu||_{H^s}.
\end{align*}
Now, we subtract two successive systems of (\ref{3-system nu}) and obtain
\begin{equation}
\begin{cases}
V^\nu = \p \Phi^\nu,\\
\lap \Phi^\nu = -\kk^{-1}A_1^{\nu-1},\\
\lap A_0^{\nu} = \kk^{-1} A_2^{\nu-1}+F_0(\p v_0^\nu)-F_0(\p v_0^{\nu-1}),\\
\lap A_1^{\nu} = \kk^{-1} A_3^{\nu-1}+F_1(\p ^{\alpha}v_0^{\nu},\p^{\beta_0}h_0^{\nu})-F_1(\p ^{\alpha}v_0^{\nu-1},\p^{\beta_0}h_0^{\nu-1}),\\
\lap A_2^{\nu} = F_2(\p^{\alpha} v_0^\nu, \p^{\beta_0}h_0^{\nu},\p^{\beta_1}h_1^{\nu})-F_2(\p^{\alpha} v_0^{\nu-1}, \p^{\beta_0}h_0^{\nu-1},\p^{\beta_1}h_1^{\nu-1}),\\
\lap A_3^{\nu} = F_3(\p^{\alpha} v_0^\nu, \p^{\beta_0}h_0^{\nu},\p^{\beta_1}h_1^{\nu},\p^{\beta_2}h_2^{\nu})- F_3(\p^{\alpha} v_0^{\nu-1}, \p^{\beta_0}h_0^{\nu-1},\p^{\beta_1}h_1^{\nu-1},\p^{\beta_2}h_2^{\nu-1}).
\label{3-system diff}
\end{cases}
\end{equation}
Here, for $k=1,2,3$, we have
\begin{dmath}
F_k(\p^{\alpha} v_0^\nu, \p^{\beta_0}h_0^{\nu},\p^{\beta_1}h_1^{\nu},\p^{\beta_2}h_2^{\nu})- F_k(\p^{\alpha} v_0^{\nu-1}, \p^{\beta_0}h_0^{\nu-1},\p^{\beta_1}h_1^{\nu-1},\p^{\beta_2}h_2^{\nu-1})\\
=C_{\alpha_1\cdots\alpha_m\beta_1\cdots\beta_n}^{\gamma_1\cdots\gamma_n,k}\Big((\p^{\alpha_1} V^\nu)\cdots(\p^{\alpha_m} v_0^\nu)(\p^{\beta_1}h_{\gamma_1}^{\nu})\cdots(\p^{\beta_n}h_{\gamma_n}^{\nu})\Big)+\cdots+\Big((\p^{\alpha_1} v_0^{\nu-1})\cdots(\p^{\alpha_m} V^\nu)(\p^{\beta_1}h_{\gamma_1}^{\nu})\cdots(\p^{\beta_n}h_{\gamma_n}^{\nu})\Big)\\
+\Big((\p^{\alpha_1} v_0^{\nu-1})\cdots(\p^{\alpha_m} v_0^{\nu-1})(\p^{\beta_1}A_{\gamma_1}^{\nu})\cdots(\p^{\beta_n}h_{\gamma_n}^{\nu})\Big)+\cdots+\Big((\p^{\alpha_1} v_0^{\nu-1})\cdots(\p^{\alpha_m} v_0^{\nu-1})(\p^{\beta_1}h_{\gamma_1}^{\nu-1})\cdots(\p^{\beta_n}A_{\gamma_n}^{\nu})\Big),
\end{dmath}
where
\begin{align*}
\alpha_1+\cdots+\alpha_m+(\beta_1+\gamma_1)+\cdots+(\beta_n+\gamma_n) = k+2,\\
1\leq \alpha_i \leq k ,\q 1\leq \beta_j+\gamma_j \leq k+1,\q 0\leq \gamma_j\leq k-1,\\
1\leq m+n\leq k+2.
\end{align*}
In addition,
\begin{dmath}
F_0(\p v_0^\nu)-F_0(\p v_0^{\nu-1}) = \delta^{ij}\delta^{kl}\Big((\p_i V_k^\nu)(\p_l v_{,0,j}^{\nu})-(\p_i v_{0,k}^{\nu-1})(\p_l V_j^{\nu})\Big).
\end{dmath}
Hence,  the same analysis which we applied to bound $m_{k}^{\nu}$ yields
\begin{align}
||V^\nu||_{H^s} = ||\p \Phi^\nu||_{H^s} \leq \kk^{-1}M_1^{\nu-1},\\
M_0^{\nu} \leq C\kk^{-1} M_2^{\nu-1}+Cq(||v_0^{\nu}||_{H^s},||v_0^{\nu-1}||_{H^s})||V^\nu||_{H^s},\\
M_1^{\nu} \leq C\kk^{-1} M_3^{\nu-1}+ Cq(m_*^{\nu},m_*^{\nu-1},||v_0^{\nu}||_{H^s},||v_0^{\nu-1}||_{H^s})(M_0^{\nu}+||V^\nu||_{H^s}),\\
M_2^{\nu} \leq Cq(m_*^{\nu},m_*^{\nu-1},||v_0^{\nu}||_{H^s},||v_0^{\nu-1}||_{H^s})(M_1^{\nu}+M_0^{\nu}+||V^\nu||_{H^s}),\\
M_3^{\nu} \leq Cq(m_*^{\nu},m_*^{\nu-1},||v_0^{\nu}||_{H^s},||v_0^{\nu-1}||_{H^s})(M_2^{\nu}+M_1^{\nu}+M_0^{\nu}+||V^\nu||_{H^s}),
\end{align}
for some polynomial $q$.

Summing these up, we have
\begin{align}
 M_*^{\nu} \leq \kk^{-1} Q(C, \mathfrak{m}, ||u_0||_{H^s}) M_*^{\nu-1},\label{sum M}
\end{align}
for some polynomial $Q$, where
$$m_*^{\nu}\leq \mathfrak{m} := \mathfrak{m}(C, ||u_0||_{H^s})$$
for each $\nu$.  Then, we have inductively that
\begin{align}
M_*^{\nu} \leq \Big(\kk^{-1} Q(C, \mathfrak{m}, ||u_0||_{H^{s-1}})\Big)^{\nu}M_*^0.
\end{align}
But since $M_*^0 = m_*^0$, and so if $\kk^{-1}$ is chosen such that
$$
\kk^{-1} Q(C, \mathfrak{m}, ||u_0||_{H^{s-1}}) < 1,
$$
then it is easy to see that
$$
M_*^{\nu}+\cdots+M_*^{\nu+n} \to 0
$$
as $\nu,n\to \infty$.

Therefore, we have proved
\thm \label{existence of data in bounded initial domain} Given the initial domain $\DD_0$ is bounded, diffeomorphic to the unit ball, and any divergence free $u_0 \in H^s, s\geq 5$, there exist data $v_0 = v_{0,\kk}$ and $h_0=h_{0,\kk}$, satisfying the compatibility condition up to order $5$, i.e.,
$$
h_k|_{\p\DD_0}=h_{k,\kk}|_{\p\DD_0}=0,\q 0\leq k\leq 5,
$$
such that the quantities
$$||v_{0,\kk}||_{H^{s}(\DD_0)} \q \text{and}\q \sum_{k=0}^r||h_{k,\kk}||_{H^{s-k}(\DD_0)},\q s\geq 5
$$ are uniformly bounded independent of $\kk$.
\rmk We give data for the enthalpy $h$ instead of the density $\rho$ in order to get bounded energy initially. If one were to do it the other way around and try to give constant $
\rho_0$ as data it would follow that $h_0$ has to be constant and hence $0$ and this would lead to that $D_t^2 h=(\partial v_0)^2$ at time $0$, and this would in general contradict that $D_t^2 h=0$ at the boundary so the compatibility conditions would not be satified and hence there would not be a solution with the required Sobolev regularity.
\rmk
Our method is systematic and so we can solve for data that satisfies $N$-compatibility conditions for any finite $N$. This, together with the fact that one can also generalize the energy estimate (\ref{main 1}) to any order, the result of \cite{L} guarantees the existence of solution for the Euler equations within $[0,T]$ for each $\kk$, which converges to the solution for the incompressible Euler equations as $\kk\to\infty$.

\subsection{Uniform bounds for $E_{4,\kk}^*(0)$}\label{section 7.3}
We are now able to show $E_{4,\kk}^*(0)$ in Section 6 is uniformly bounded regardless of $\kk$. This is because
\begin{align}
\sum_{k+s\leq 4}\int_{\Omega}\rho_0Q(\p^sh_k,\p^sh_k)\dx \lesssim \sum_{0\leq k\leq 4}||h_k||_{H^{4-k}(\Omega)}^2 \leq \mathfrak{m},
\end{align}
and by the trace lemma,
\begin{dmath}
\sum_{k+s\leq 4}\int_{\p\Omega}\rho_0Q(\p^sh_k,\p^sh_k)\,dS \lesssim \sum_{0\leq k\leq 4}||h_k||_{H^{4-k}(\p\Omega)}^2 \lesssim  \sum_{0\leq k\leq 4}||h_k||_{H^{5-k}(\Omega)}^2 \leq \mathfrak{m}.
\end{dmath}
In addition to these, we have
\begin{dmath}
\sum_{0\leq k\leq 4}||(\p^{4-k}D_t^k v)\big|_{t=0}||_{L^2} \lesssim ||v_0||_{H^4}+\sum_{0\leq k\leq 3}||h_k||_{H^{4-k}},
\end{dmath}
since $D_t v = -\p h$. This shows
\begin{align}
\sum_{k+s\leq 4}\int_{\Omega}\rho_0Q(\p^s D_t^k v\big|_{t=0},\p^s D_t^kv\big|_{t=0})\dx
\end{align}
is uniformly bounded as well. Finally, since $h_4 = h_5=0$ in $\Omega$, we have
\begin{align*}
\sum_{k=2}^5 W_k(0) \lesssim \mathfrak{m},
\end{align*}
 and hence we have $E_{4,\kk}^*(0)$ bounded uniformly.
\subsection{The convergence of $v_{0,\kk}$ and $h_{0,\kk}$ and the physical sign condition}
Since for each $s\geq 5$, we have the estimate
\begin{align}
||v_{0,\kk}-u_0||_{H^s} \leq ||\p \phi_{\kk}||_{H^s} \lesssim \kk^{-1}||h_{1,\kk}||_{H^{s-1}}.
\end{align}
Because of this, we have $||v_{0,\kk}-u_0||_{C^1}\leq ||v_{0,\kk}-u_0||_{H^s}\lesssim \kk^{-1}||h_{1,\kk}||_{H^{s-1}}$ whenever $s>\frac{n}{2}+1$, which implies $v_{0,\kk}\to u_0$ in $C^1$ since $||h_{1,\kk}||_{H^s}$ is bounded uniformly independent of $\kk$.
On the other hand, since $\DD_0$ is bounded, we assume the physical sign condition holds when $t=0$, i.e.,
\begin{align}
\nab_N h_0 \leq -\epsilon <0,\q \text{on}\,\,\, \p\DD_0. \label{taylor sign at 0}
\end{align}
 This will be true under small perturbation in $[0,T]$ due to \eqref{check E}.
 Given any data for the incompressible equations $u_0$ such that the corresponding $p_0$ satisfies $-\nab_N p_0\geq \epsilon>0$, our data for the compressible equations $h_{0,\kk}$ will also satisfy \eqref{taylor sign at 0} if
 $\kk^{-1}$ is sufficiently small. In fact, since
 $$\triangle h_{0,\kk}=(\partial v_{0,\kk})^2+\kk^{-1} h_{2,\kk},$$
and so
$$\triangle (h_{0,\kk}-p_0)=\kk^{-1} h_{2,\kk}+(\p^2 \phi_{\kk})(\p u_0)+(\p^2 \phi_{\kk})^2.$$
On the other hand, the standard elliptic estimates yield
$$\|h_{0,\kk}-p_0\|_{H^s}\lesssim \kk^{-1} \|h_{2,\kk}\|_{H^{s-2}}+\kk^{-1}||u_0||_{H^s}||h_1||_{H^{s-2}}.$$
This yields the convergence of $h_{0,\kk}\to p_0$ in $C^1$,and so Theorem \ref{convergence of data} is proved.

\begin{appendix}
\section{Appendix}
\noindent {\bf List of notations:}

\begin{itemize}
\item $D_{t}$: the material derivative
\item $\p_i$: partial derivative with respect to Eulerian coordinate $x_i$
\item $\DD_t\in\RR^n$: the domain occupied by fluid particles at time $t$ in Eulerian coordinate
\item $\Omega\in\RR^n$: the domain occupied by fluid particles in Lagrangian coordinate
\item $\p_a = \frac{\p}{\p y_a}$: partial derivative with respect to Lagrangian coordinate $y_a$
\item $\nab_a$: covariant derivative with respect to $y_a$
\item $\Pi S$: projected tensor $S$ on the boundary
\item $\cnab,\cp$: projected derivative on the boundary
\item $N$: the outward unit normal of the boundary
\item $\theta=\cnab N$: the second fundamental form of the boundary
\item $\sigma = tr(\theta)$: the mean curvature
\end{itemize}
{\bf Mixed norms}
\begin{itemize}
\item $\lee\cdot\ree_r = \sum_{k+s=r}||\nab^sD_t^k\cdot||_{\llb}$
\item $||\cdot||_{r,0} = \sum_{s+k=r,k<r}||\nab^s D_t^k\cdot||_{\lli}$
\item $||\cdot||_{r}= ||\cdot||_{r,0}+||\sqrt{e'(h)}D_t^r\cdot||_{\lli}$

\end{itemize}

\subsection{Extension of the normal to the interior and the geodesic normal coordinate}
\indent The definition of our energy (\ref{Er}) relies on extending the normal to the interior, which is done by foliating the domain close to the boundary into the surface that do not self-intersect.  We also want to control the time evolution of the boundary, which can be measured by the time derivative of the normal in the Lagrangian coordinate. We conclude the above statements by the following two lemmas, whose proof can be found in \cite{CL}.\\
\lem \label{trace 1} let $l_0$ be the injective radius (\ref{inj rad}), and let $d(y)=dist_g(y,\p\Omega)$ be the geodesic distance in the metric $g$ from $y$ to $\p\Omega$. Then the co-normal $n=\nab d$ to the set $S_a=\p\{y\in\Omega:d(y)=a\}$ satisfies, when $d(y)\leq \frac{l_0}{2}$ that
\begin{align}
|\nab n| \lesssim |\theta|_{L^{\infty}(\p\Omega)}, \label{extending nab n}\\
|D_t n| \lesssim |D_t g|_{L^{\infty}(\Omega)},
\end{align}
where we have used the convention that $A\lesssim B$ means $A\leq CB$ for universal constant $C$.
\lem \label{trace 2} let $l_0$ be the injective radius (\ref{inj rad}),and let $d_0$ be a fixed number such that $\frac{l_0}{16}\leq d_0\leq \frac{l_0}{2}$. Let $\eta$ be a smooth cut-off function satisfying $0\leq \eta(d)\leq1$, $\eta(d)=1$ when $d\leq\frac{d_0}{4}$ and $\eta(d)=0$ when $d>\frac{d_0}{2}$. Then the psudo-Riemannian metric $\gamma$ given by
$$
\gamma_{ab}=g_{ab}-\tilde{n}_a\tilde{n}_b,
$$
where $\tilde{n}_c=\eta(\frac{d}{d_0})\nab_cd$
satisfies
\begin{align}
|\nab\gamma|_{L^{\infty}(\Omega)}\lesssim(|\theta|_{L^{\infty}(\p\Omega)}+\frac{1}{l_0})\\
|D_t\gamma(t,y)|\lesssim |D_t g|_{L^{\infty}(\Omega)}. \label{Dtgamma}
\end{align}
\rmk The above two lemmas yield that the quantities $|D_t n|$ and $|D_t\gamma(t,y)|$ involved in the $Q$-inner product is controlled by the a priori assumptions (\ref{geometry_bound})-(\ref{Dtth}),since $D_tg$ behaves like $\nab v$ by (\ref{Dtg}). Hence, the time derivative on the coefficients of the $Q$-inner product generates only lower-order terms. In addition, by (\ref{geometry_bound}) , $|\nab n|$ and $|\nab\gamma|$ are controlled by $K$, which is essential when proving the elliptic estimates.

\subsection{Sobolev lemmas}
Let us now state some Sobolev lemmas in a domain with boundary.
\lem \label{interior sobolev} (Interior Sobolev inequalities)
Suppose $\frac{1}{l_0}\leq K$ and $\alpha$ is a $(0,r)$ tensor, then
\begin{align}
||\alpha||_{L^{\frac{2n}{n-2s}}(\Omega)}\lesssim_{K} \sum_{l=0}^s||\nab^l\alpha||_{\lli},\q 2s<n \label{interior sobolev 2s<n},\\
||\alpha||_{L^{\infty}(\Omega)}\lesssim_{K} \sum_{l=0}^{s}||\nab^l\alpha||_{\lli},\q 2s>n \label{interior sobolev 2s>n}.
\end{align}
Similarly, on $\p\Omega$, we have
\lem \label{boundary soboolev} (Boundary Sobolev inequalities)
\begin{align}
||\alpha||_{L^{\frac{2(n-1)}{n-1-2s}}(\Omega)}\lesssim_{K} \sum_{l=0}^s||\nab^l\alpha||_{\llb},\q 2s<n-1,\\
||\alpha||_{L^{\infty}(\Omega)}\lesssim_{K} \delta||\nab^s\alpha||_{\llb}+\delta^{-1}\sum_{l=0}^{s-1}||\nab^l\alpha||_{\llb},\q 2s>n-1,
\end{align}
for any $\delta>0$. In addition, for the boundary we can also interpret the norm be given by the inner product $\langle \alpha, \alpha\rangle=\gamma^{IJ}\alpha_I\alpha_J$, and the covariant derivative is then given by $\cnab$.
\subsection{Interpolation on spatial derivatives}
We shall first record spatial interpolation inequalities. Most of the results are are standard in $\RR^n$, but we must control how it depends on the geometry of our evolving domain. The coefficients involved in our inequalities depend on $K$, whose reciprocal is the lower bound for the injective radius $l_0$. We omit the proofs, which can be found in the appendix of \cite{CL}.
\thm \label{int interpolation} (Interior interpolation)
Let $u$ be a $(0,r)$ tensor, and suppose $\frac{1}{l_0}\leq K$, we have
\begin{align}
\sum_{j=0}^{l}||\nab^ju||_{L^{\frac{2r}{k}}(\Omega)}\lesssim ||u||_{L^{\frac{2(r-l)}{k-l}}(\Omega)}^{1-\frac{l}{r}}(\sum_{i=0}^{r}||\nab^iu||_{\lli}K^{r-i})^{\frac{l}{r}}.
\end{align}
In particular, if $k=l$,
\begin{align}
\sum_{j=0}^{k}||\nab^ju||_{L^{\frac{2r}{k}}(\Omega)}\lesssim ||u||_{L^{\infty}(\Omega)}^{1-\frac{k}{r}}(\sum_{i=0}^{r}||\nab^iu||_{\lli}K^{r-i})^{\frac{k}{r}}. \label{31}
\end{align}

\subsection{Interpolation on $\p\Omega$}
We need the following boundary interpolation inequalities to estimate the boundary part of our energy (\ref{Er}).
\thm \label{bdy interpolation} (Boundary interpolation)
Let $u$ be a $(0,r)$ tensor, then
\begin{equation}
||\cnab^l u||_{L^{\frac{2r}{k}}(\p\Omega)}\lesssim ||u||_{L^{\frac{2(r-l)}{k-l}}(\p\Omega)}^{1-\frac{l}{r}}||\cnab^ru||_{\llb}^{\frac{l}{r}}.
\end{equation}
In particular, if $k=l$,
\begin{align}
||\cnab^k u||_{L^{\frac{2r}{k}}(\p\Omega)}\lesssim ||u||_{L^{\infty}(\p\Omega)}^{1-\frac{k}{r}}||\cnab^ru||_{\llb}^{\frac{k}{r}}.
\end{align}

\thm \label{gag-ni thm}(Gagliardo-Nirenberg interpolation inequality)
Let $u$ be a $(0,r)$ tensor, and suppose $\p\Omega\in\RR^2$ and $\frac{1}{l_0}\leq K$, we have
\begin{align}
||u||_{L^4{(\p\Omega})}^2 \lesssim_K ||u||_{\llb}||u||_{H^1(\p\Omega)},\label{gag-ni}
\end{align}
where the boundary Sobolev norm $||u||_{H^1(\p\Omega)}$ is defined via tangential derivative $\cnab$.
\begin{proof}
It suffices for us to work in the local coordinate charts $(U_i)_{i=1}^n$ of $\p\Omega$. We introduce the corresponding partition of unity $(\chi_i)_{i=1}^n$, where each $\chi_i$ is supported in $U_i$ and vanishing on the boundary of $U_i$. Now by the result of Constantin and Seregin \cite{Cons}, we have
$$
||u_i||_{L^4(U_i)}^2 \lesssim ||u_i||_{L^2(U_i)}||\cnab u_i||_{L^2(U_i)},
$$
where $u_i = \chi_iu$. But since
\begin{equation}
||\cnab u_i||_{L^2(U_i)} = ||\cnab (\chi_i u)||_{L^2(U_i)} \leq |\cnab \chi_i|_{L^{\infty}}||u||_{L^2(U_i)}+||\chi_i\cnab u||_{L^2(U_i)},\label{gag local}
\end{equation}
(\ref{gag-ni}) follows by summing up (\ref{gag local}) since $\chi_i$ can be chosen so that
$
\sum_i|\nab \chi_i| \leq C(K),
$
as long as $\frac{1}{l_0}\leq K$ (see \cite{CL}).
\end{proof}
\rmk One can also prove a generalized (\ref{gag-ni}) of the form
\begin{equation}
||u||_{L^p{(\p\Omega})}^2 \lesssim ||u||_{L^{p/2}(\p\Omega)}||u||_{H^1(\p\Omega)},\qquad p\geq 4.
\end{equation}

\indent Our next theorem shall be dealing with the projected derivatives acting on tensors. We first define that if $\alpha$ is a $(0,t)$ tensor, then the projected $(0,r), r<t$ derivative $\Pi^{r,0}\nab^r \alpha$ has components
$$(\Pi\nab^r)_{i_1,\cdots,i_r} \alpha_{i_r+1,\cdots,i_t}=\gamma_{i_1}^{j_1}\cdots\gamma_{i_r}^{j_r}\nab_{j_1}\cdots\nab_{j_r} \alpha_{i_r+1,\cdots,i_t}.$$
The proof of the next theorem is rather involved, we refer \cite{CL} for the full proof.
\thm \label{tensor interpolation} (Tensor interpolation)
Let $\alpha$ be a $(0,t)$ tensor, and let $r'= r-2$, suppose that $|\theta|+|\frac{1}{l_0}|\leq K$, then for $t+s<r$ we have
\begin{multline}
||(\Pi^{s,0}\nab^s)\alpha||_{L^{\frac{2r'}{s}}(\p\Omega)}\lesssim_K ||\alpha||_{\linf}^{1-s/r'}\Big(||\nab^{r'}\alpha||_{\llb}+(1+||\theta||_{\linf})^{r'}\\
\cdot(||\theta||_{\linf}+||\cnab^{r'}\theta||_{\llb})\sum_{l=0}^{r'-1}||\nab^l\alpha||_{\llb}\Big)\\
+(1+||\theta||_{\linf})^s(||\theta||_{\linf}+||\cnab^{r'}\theta||_{\llb})^{s/r'}\sum_{l=0}^{r'-1}||\nab^l\alpha||_{\llb}.
\end{multline}
In particular,
\begin{multline}
\Big{\|}|(\Pi^{s,0}\nab^s)\alpha||(\Pi^{r'-s,0}\nab^{r'-s}\beta)|\Big{\|}_{\llb}\lesssim_K(||\alpha||_{\linf}+\sum_{l=0}^{r'-1}||\nab^l\alpha||_{\llb})||\nab^{r'}\beta||_{\llb}\\
+(||\beta||_{\linf}+\sum_{l=0}^{r'-1}||\nab^l\beta||_{\llb})||\nab^{r'}\alpha||_{\llb}+(1+||\theta||_{\linf})^{r'}(||\theta||_{\linf}+||\cnab^{r'}\theta||_{\llb})\\
+(||\alpha||_{\linf}+\sum_{l=0}^{r'-1}||\nab^l\alpha||_{\llb})(||\beta||_{\linf}+\sum_{l=0}^{r'-1}||\nab^l\beta||_{\llb}) \label{tensor int}.
\end{multline}
\begin{proof}
See \cite{CL}, section 4.
\end{proof}
\subsection{The trace theorem}
\thm \label{trace theorem} (Trace theorem)
Let $\alpha$ be a $(0,r)$ tensor, and assume that $|\theta|_{\linf}+\frac{1}{l_0}\!\leq \!K$. Then
\begin{align}
||\alpha||_{\llb}\lesssim_{K,r,n} \sum_{j\leq 1}||\nab^j\alpha||_{\lli} \label{trace}
\end{align}
\begin{proof}
Let $\mathcal{N}$ be the extension of the normal to the interior, then the Green's identity yields
$$
\int_{\p\Omega}|\alpha|^2\,d\mu_{\gamma} = \int_{\Omega}\nab_k(\mathcal{N}^k|\alpha|^2)\,d\mu.
$$
Hence, by Lemma \ref{trace 1} and \ref{trace 2}, (\ref{trace}) follows.
\end{proof}

\end{appendix}

\section*{Acknowledgments}
H.L \ is partially supported by NSF grant DMS--1237212. The authors would like to thank the referee for his/her fruitful comments on the previous version of this manuscript.

\end{document}